\documentclass[12pt]{amsart}
\usepackage{latexsym}
\usepackage{amsthm}
\usepackage{amsmath}
\usepackage[dvips]{graphics}
\usepackage{graphicx}
\usepackage{amssymb}
\usepackage{color}
\usepackage{ulem}
\usepackage{epsfig}
\usepackage[all]{xy}

\begin{document}

\numberwithin{equation}{section}

\theoremstyle{plain}
\newtheorem{theorem}{Theorem}[section]
\newtheorem{lemma}[theorem]{Lemma}
\newtheorem{conj}[theorem]{Conjecture}
\newtheorem{pposition}[theorem]{Proposition}
\newtheorem{remark}[theorem]{Remark}

\newtheorem*{theorem*}{Theorem}
\newtheorem*{lemma*}{Lemma}
\newtheorem*{claim*}{Claim}
\newtheorem*{exercise*}{Exercise}
\newtheorem*{note*}{Note}
\newtheorem*{example*}{Example}
\newtheorem*{problem*}{Problem}
\newtheorem*{solution*}{Solution}
\newtheorem*{remark*}{Remark}

\newtheorem{corollary}[theorem]{Corollary}
\newtheorem{example}[]{Example}
\newtheorem{conjecture}[theorem]{Conjecture}
\newtheorem{definition}[theorem]{Definition}
\newtheorem{proposition}[theorem]{Proposition}

\newcommand{\nature}{\ensuremath{\mathbb{N}}}   % natural numbers
\newcommand{\integer}{\ensuremath{\mathbb{Z}}}  % integer
\newcommand{\cross}{\ensuremath{^{\times}}}     % cross
\newcommand{\rational}{\ensuremath{\mathbb{Q}}} % rational
\newcommand{\Z}{\integer}
\newcommand{\f}{\frac}
\newcommand{\eq}{\equiv}
\newcommand{\N}{\mathbb N}
\newcommand{\mo}{{\rm mod}}
\newcommand{\zpinfty}{\integer (p^\infty)}
\newcommand{\x}{\text{\boldmath{$X$}}}
\newcommand{\h}{\mathbb{H}}
\newcommand{\HH}{\text{\boldmath{$H$}}}
\newcommand{\g}{\text{\boldmath{$g$}}}
\newcommand{\G}{\text{\boldmath{$G$}}}
\newcommand{\R}{\text{ $R$ }}
\newcommand{\SF}{\mathcal{SF}}
\newcommand{\F}{\text{\boldmath{$F$}}}
\newcommand{\w}{\text{\boldmath{$w$}}}
\newcommand{\re}{\text{Re}}
\newcommand{\im}{\text{Im}}
\newcommand{\vv}{\text{\boldmath{$v$}}}
\newcommand{\order}{\text{\boldmath{$O$}}}
\newcommand{\iso}{\cong}
\renewcommand{\thefootnote}{}

\hbox{Trans. Amer. Math. Soc. 362(2010), no.\,12, 6425--6455.}
\medskip
\title{On Almost Universal Mixed Sums of Squares and Triangular Numbers}
\author{Ben Kane}
\address{Department of Mathematics, University of Cologne, Weyertal 86-90, 50931 Cologne, Germany}
\email{bkane@math.uni-koeln.de}
\author{Zhi-Wei Sun}
\address{Department of Mathematics\\
Nanjing University\\
Nanjing 210093\\
People's Republic of China}
\email{zwsun@nju.edu.cn}
\thanks {On August 20, 2008, the initial version was posted to
{\tt http://arxiv.org/abs/0808.2761}.
\\\indent
The second author was the corresponding author, and he was supported
by the National Natural Science Foundation (grant 10871087) of P. R.
China.  This research was conducted when the first author was a
postdoctor at Radboud Universiteit, Nijmegen, Netherlands.}

\begin{abstract} In 1997 K. Ono and K. Soundararajan [Invent. Math. 130(1997)]
proved that under the generalized Riemann hypothesis
any positive odd integer greater than $2719$
can be represented by the famous Ramanujan form $x^2+y^2+10z^2$,
equivalently the form $2x^2+5y^2+4T_z$ represents all integers greater than 1359,
where $T_z$ denotes the triangular number $z(z+1)/2$.
Given positive integers $a,b,c$ we employ modular
forms and the theory of quadratic forms to determine completely
when the general form
$ax^2+by^2+cT_z$ represents
sufficiently large integers and establish similar results for
the forms $ax^2+bT_y+cT_z$ and $aT_x+bT_y+cT_z$.
Here are some consequences of our main theorems:
(i) All sufficiently large odd numbers have the form $2ax^2+y^2+z^2$
if and only if all prime divisors of $a$ are congruent to 1 modulo 4.
(ii) The form
$ax^2+y^2+T_z$ is almost universal (i.e., it represents sufficiently
large integers) if and only if each odd prime divisor of $a$ is
congruent to 1 or 3 modulo 8.
(iii) $ax^2+T_y+T_z$
is almost universal if and only if all odd prime divisors of $a$ are
congruent to 1 modulo 4.
(iv) When $v_2(a)\not=3$, the form
$aT_x+T_y+T_z$ is almost universal if and only if all odd prime
divisors of $a$ are congruent to 1 modulo 4 and
$v_2(a)\not=5,7,\ldots$, where $v_2(a)$ is the $2$-adic order of
$a$.
\end{abstract}
 \maketitle

\section{Introduction and the Main Results}\label{introsection}

A classical theorem of Lagrange states that any
$n\in\N=\{0,1,2,\ldots\}$ can be written as a sum of four squares (of integers).
In 1916 S. Ramanujan \cite{Ramanujan} found all the finitely many vectors $(a,b,c,d)$ with
$a,b,c,d\in\Z^+=\{1,2,3,\ldots\}$ such that the form $ax^2+by^2+cz^2+dw^2$ (with $x,y,z,w\in\Z$)
represents all natural numbers. Ramanujan also asked for determining those vectors
$(a,b,c,d)\in(\Z^+)^4$ such that the form $ax^2+by^2+cz^2+dw^2$ represents all sufficiently large
integers; this problem was essentially solved by H. D. Kloosterman \cite{Kloo} with help from the useful
Kloosterman sum, and this work represents a major breakthrough in the field
of quadratic forms.

What about ternary quadratic forms? A well known theorem of Gauss and Legendre states that
$n\in\N$ is a sum of three squares if and only if it is not of the form $4^k(8l+7)$ with $k,l\in\N$.
In general, it is known that for any $a,b,c\in\Z^+$ the subset $\{ax^2+by^2+cz^2:\ x,y,z\in\Z\}$
of $\N$ cannot have asymptotic density 1 because there is always a congruence class
modulo a power of some prime $p$ dividing $2abc$
which is not even locally represented by the form $ax^2+by^2+cz^2$.

For $x\in\Z$ let $T_x$ denote the triangular number $x(x+1)/2$.
Clearly $T_n=T_{-n-1}$ for all $n\in\N$. A famous
assertion of Fermat states that each $n\in\N$ can be
expressed as a sum of three triangular numbers,
equivalently $8n+3$ is a sum of three (odd) squares;
this follows immediately from the Gauss-Legendre theorem.
Here is another consequence of the Gauss-Legendre theorem
observed by Euler: Each natural number can be written in the
form $x^2+y^2+T_z$ with $x,y,z\in\Z$. Recently, B. K. Oh and the second author \cite{Sun2}
showed that for any $n\in\Z^+$ there are $x,y,z\in\Z$
such that $n=x^2+(2y+1)^2+T_z$, i.e., $n-1=x^2+8T_y+T_z$.

In view of the above, it is natural to study mixed sums of squares and triangular
numbers of the following three types:
$$ax^2+by^2+cT_z,\ ax^2+bT_y+cT_z,\ aT_x+bT_y+cT_z$$
where $a,b,c\in\Z^+=\{1,2,3,\ldots\}$. Let $f(x,y,z)$ be any of
the three forms, and define the exceptional set
$$E(f):=\{n\in\N: f(x,y,z)=n\text{ has no integral solution}\}.$$
If $E(f)=\emptyset$, then $f$ is said to be {\it universal}; if
$E(f)$ is finite, then we call $f$ {\it almost universal}. When the
set $E(f)$ has asymptotic density zero, i.e.,
$$\lim_{N\to+\infty}\frac{|\{1\le n\le N:\ f(x,y,z)=n\ \text{for some}\ x,y,z\in\Z\}|}N=1,$$
we say that $f$ is {\it asymptotically universal}. In the case $\gcd(a,b,c)>1$,
obviously $f$ is neither almost universal nor asymptotically universal.

In 1862 J. Liouville (cf. \cite[p.\,23]{Dickson}) proved the following
result: For positive integers $a\le b\le c$, the form
$aT_x+bT_y+cT_z$ is universal if and only if $(a,b,c)$ is among
the following vectors:
$$(1,1,1),\ (1,1,2),\ \ (1,1,4),\ (1,1,5),\ (1,2,2),
\ (1,2,3),\ (1,2,4).$$ Recently the second author \cite{Sun1}
initiated the determination of all universal forms of the type $ax^2+by^2+cT_z$ or
$ax^2+bT_y+cT_z$, and the project was finally completed by combining the results in
\cite{Sun1}, \cite{GuoSun1} and \cite{Sun2}. Namely, for $a,b,c\in\Z^+$ with $a\le b$,
the form $ax^2+by^2+cT_z$ is universal if and only if $(a,b,c)$ is among
the following vectors:
\begin{align*}(1,1,1),\ (1,1,2),\ (1,2,1),\
(1,2,2),\ (1,2,4),
\\(1,3,1),\ (1,4,1),\ (1,4,2),\ (1,8,1),\ (2,2,1).
\end{align*}
Also, for $a,b,c\in\Z^+$ with  $b\ge c$, the form $ax^2+bT_y+cT_z$ is
universal if and only if $(a,b,c)$ is among the following vectors:
\begin{align*}&(1,1,1),\ (1,2,1),\ (1,2,2),\ (1,3,1),\ (1,4,1),\
(1,4,2),\ (1,5,2),
\\ &(1,6,1),\ (1,8,1),
\ (2,1,1),\ (2,2,1),\ (2,4,1),\ (3,2,1), \ (4,1,1),\ (4,2,1).
\end{align*}

In 1916 Ramanujan (cf. \cite{Ramanujan} and \cite{Ono2}) conjectured that those positive even integers
not represented by $x^2+y^2+10z^2$
are exactly those of the form $4^k(16l+6)$ with $k,l\in\N$, and that
those positive odd integers not represented by $x^2+y^2+10z^2$ are as follows:
$$3,\,7,\,21,\,31,\,33,\,43,\,67,\,79,\,87,\,133,\,217,\,219,\,223,\,253,\,307,\,391,\ldots.$$
In 1927 L. E. Dickson \cite{Dickson27} proved Ramanujan's conjecture about even numbers
by a simple argument. However, Ramanujan's conjecture about odd numbers is very difficult.
For $n\in\N$, clearly
\begin{align*} &2n+1=x^2+y^2+10z^2\ \text{for some}\ x,y,z\in\Z
\\\iff &2n+1=(2x)^2+10y^2+(2z+1)^2\ \text{for some}\ x,y,z\in\Z
\\\iff&n=2x^2+5y^2+4T_z\ \text{for some}\ x,y,z\in\Z.
\end{align*}
Only in 1990 were W. Duke and R. Schulze-Pillot \cite{DSP} able to
show that sufficiently large odd integers can be written in the form
$x^2+y^2+10z^2$, or equivalently that the form $2x^2+5y^2+4T_z$ is
almost universal. In 1997 K. Ono and K. Soundararajan
\cite{OnoSound} showed further that the generalized Riemann
hypothesis implies that the only positive odd integers not in the
form $x^2+y^2+10z^2$ are those listed by Ramanujan together with
$679$ and $2719$, in other words $E(2x^2+5y^2+4T_z)$ consists of the
following numbers:
$$1,\,3,\,10,\,15,\,16,\,21,\,33,\,39,\,43,\,66,
\,108,\,109,\,111,\,126,\,153,\,195,\,339,\,1359.$$

Motivated by his conjecture on sums of primes and triangular numbers
(cf. \cite[Conjecture 1.1]{Sun4}),
the second author \cite{Sun3} recently conjectured
that for any $k,l\in\N$ the form $2^kx^2+2^ly^2+T_z$ is almost universal.
The first author \cite{Kane4} showed that all of those forms conjectured to be almost universal in \cite{Sun3}
are asymptotically universal and many of them are
almost universal.

In this paper we aim at determining all
asymptotically universal forms and almost universal forms of the three types
via modular forms and the theory of quadratic forms.

For convenience we introduce some basic notation. We may write a
positive integer $a$ in the form $2^{v_2(a)}a'$ with $v_2(a)\in\N$
and $a'$ odd; $v_2(a)$ is called the $2$-adic order of $a$
(equivalently, $2^{v_2(a)}\|a$) while $a'$ is said to be the odd
part of $a$. For $a\in\Z$ and $m\in\Z^+$, by $a\R m$ we mean that
$a$ is quadratic residue modulo $m$, i.e., $a$ is relatively prime
to $m$ and the equation $x^2\equiv a\ (\mo\ m)$ is solvable over
$\Z$. For an integer $a$ and a positive odd integer $m$, it is well known that
$a\R m$ if and only if the Legendre symbol $(\f ap)$ equals 1 for every prime divisor $p$ of $m$.

Now we state our results on asymptotically universal forms.

\begin{theorem}\label{twoalmost} Fix $a,b,c\in \integer^{+}$ with $\gcd(a,b,c)=1$. Then the form
$$f(x,y,z):=ax^2+by^2 + cT_z$$
is asymptotically universal if and only if we have the following $(1)-(2)$.
\begin{enumerate}
\item $ -2bc \R a',\ -2ac \R b'\text{, and }-ab\R c'. $
\item
Either $4\nmid c$, or both $4\| c$ and $2\| ab$.
\end{enumerate}
\end{theorem}

\begin{theorem}\label{onealmost} Fix  $a,b,c\in \integer^{+}$ with $\gcd(a,b,c)=1$. Then the form
$$
f(x,y,z):=ax^2+b T_y + cT_z
$$
is asymptotically universal if and only if we have the following
$(1)-(2)$.
\begin{enumerate}
\item  $ -bc \R a',\ -2ac \R b'\text{, and }-2ab\R c'. $
\item Either $4\nmid b$ or $4\nmid c$.
\end{enumerate}
\end{theorem}

\begin{theorem}\label{trialmost}
Fix $a,b,c\in \integer^{+}$ with $\gcd(a,b,c)=1$. Then the form
$$
f(x,y,z):=aT_x+b T_y + cT_z
$$
is asymptotically universal if and only if
$$-bc \R a',\ -ac \R b'\text{, and }-ab\R c'.$$
\end{theorem}

\noindent{\it Remark}\ 1.1. For $a,b,c\in\Z^+$, if the form $ax^2+by^2+cT_z$ or $ax^2+bT_y+cT_z$ or $aT_x+bT_y+cT_z$
is asymptotically universal then $a',b',c'$ must be pairwise coprime by Theorems 1.1-1.3.

\medskip

The law of quadratic reciprocity gives restrictions under which the relations in the above theorems cannot occur.
\begin{corollary}\label{quadreccor}
Fix $a,b,c\in \integer^{+}$ and consider
\begin{align*}
(1)\ ax^2+by^2+cT_z,\qquad (2)\ aT_x +bT_y+cz^2, \qquad (3)\ aT_x + bT_y+cT_z,\\
(4)\ ax^2+bT_y+cz^2,\qquad (5)\ ax^2 +bT_y+cT_z.
\end{align*}

{\rm (i)} If $a'\eq b'\eq -c'\pmod{8}$, then none of $(1)$-$(5)$ is asymptotically universal.

{\rm (ii)} If
$$\begin{cases}a'\eq b'\eq c'+4\pmod{8}&\\v_2(a)\not\eq v_2(b)\pmod{2}\end{cases}
\quad\text{or}\quad
\begin{cases}\pm a'\eq -b'\eq c'+4\pmod{8}&\\v_2(a)\equiv v_2(b)\pmod{2},\end{cases}$$
then none of $(1)$-$(3)$ is asymptotically universal.

{\rm (iii)} If
$$\begin{cases}a'\eq b'\eq c'+4\pmod{8}&\\v_2(a)\eq v_2(b)\pmod{2}
\end{cases}\quad\text{or}\quad
\begin{cases}\pm a'\eq -b'\eq c'+4\pmod{8}&\\v_2(a)\not\eq v_2(b)\pmod{2},
\end{cases}$$
 then neither $(4)$ nor $(5)$ is asymptotically universal.
\end{corollary}

\begin{corollary}\label{new} Let $a,b,c\in\Z^+$ with $v_2(b)\eq v_2(c)\ (\mo\ 2)$.
Assume that $a'\eq b'\eq c'\ (\mo\ 4)$ fails. Then,
either none of the forms $ax^2+by^2+2cT_z$ and $ax^2+cy^2+2bT_z$ is asymptotically universal,
or none of the forms $ax^2+2cy^2+bT_z$ and $ax^2+2by^2+cT_z$ is asymptotically universal.
\end{corollary}

Any $n\in\Z^+$ can be uniquely written in the form $a^2q$ with $a,q\in\Z^+$ and $q$ squarefree, and we use
$\SF(n)$ to denote $q=\prod_{p\mid n,\,2\nmid v_p(n)}p$, the squarefree part of $n$.

Now we turn to almost universal forms.

\begin{theorem}\label{twosufflarge}
Let $a,b,c\in \integer^{+}$ with $\gcd(a,b,c)=1$ and $v_2(a)\ge v_2(b)$.
Suppose that both $(1)$ and $(2)$ in Theorem \ref{twoalmost} hold.
Then there are infinitely many positive integers not represented by the form
$$
f(x,y,z):=ax^2+by^2 + cT_z
$$
(i.e., $f$ is not almost universal) if and only if we have the following $(1)-(3)$.
\begin{enumerate}
\item \label{abprimecond} $2|a$, $4\nmid c$, $a'\equiv b'\pmod{2^{3-v_2(c)}}$, and
$$\begin{cases}4\nmid b\Rightarrow v_2(a)\eq c\ (\mo\ 2)&
\\2\nmid bc\Rightarrow 8\mid a\ \&\ 8\mid(b-c).
\end{cases}$$

\item \label{primeprodcond} All prime divisors
of $\SF(a'b'c')$ are congruent to $1$ modulo $4$ if $v_2(a)\eq v_2(b)\ (\mo\ 2)$,
and congruent to $1$ or $3$ modulo $8$ otherwise.
\item \label{repcond} $2^{3-v_2(c)}(ax^2+by^2) +c'z^2 =\SF(a'b'c')$ has no integral solutions.
\end{enumerate}
\end{theorem}

\medskip
\noindent{\it Remark} 1.2. When $ax^2+by^2+cT_z$ (with $a,b,c\in\Z^+$) is asymptotically universal
it is not necessary that (2) in Theorem \ref{twosufflarge} holds.
For example, $6x^2+y^2+10T_z$ is asymptotically universal
but we don't have (2) in Theorem \ref{twosufflarge} with $a=6$, $b=1$ and $c=10$.

\medskip

{\it Example} 1.1. Consider those forms $ax^2+by^2+cT_z$ with
$a,b,c\in\Z^+$ and $a+b+c\le10$. By Theorem \ref{twosufflarge}, we
find that those asymptotically universal ones are all almost universal.
Below is a complete list of those forms $ax^2+by^2+cT_z$ with
$a,b,c\in\Z^+$ and $a+b+c\le10$ which are almost universal but not
universal:
\begin{align*}&x^2+2y^2+3T_z,\ 2x^2+4y^2+T_z,\  x^2+5y^2+2T_z,\  x^2+6y^2+T_z,
\\&x^2+y^2+5T_z,\ 2x^2+3y^2+2T_z,\ \ 2x^2+5y^2+T_z,\ 3x^2+4y^2+T_z,
\\&x^2+2y^2+6T_z,\ x^2+5y^2+3T_z,\ 2x^2+2y^2+5T_z,\ 2x^2+4y^2+3T_z,
\\&4x^2+4y^2+T_z,\ x^2+4y^2+5T_z,\ 2x^2+3y^2+5T_z.
\end{align*}
For the four forms in the first row, the second author \cite{Sun1}
conjectured that
\begin{align*}&E(x^2+2y^2+3T_z)=\{23\},\ E(2x^2+4y^2+T_z)=\{20\},
\\&E(x^2+5y^2+2T_z)=\{19\},\ E(x^2+6y^2+T_z)=\{47\},
\end{align*}
which was confirmed by the first author \cite{Kane4} under the
generalized Riemann hypothesis.  For the form $4x^2+4y^2+T_z$ the
second author \cite{Sun4} conjectured that $E(4x^2+4y^2+T_z)$
consists of the following 19 numbers:
\begin{align*} &2,\ 12,\ 13,\, 24,\ 27,\ 34,\ 54,\ 84,\ 112,\ 133,
\\&162,\ 234,\ 237,\ 279,\ 342,\  399,\ 652,\  834,\ 864.\end{align*}
For the ten remaining forms on the above list, our computation via
computer suggests the following information:
\begin{align*}&E(x^2+y^2+5T_z)=\{3,\, 11,\, 12,\, 27,\, 129,\, 138,\, 273\},
\\&E(2x^2+3y^2+2T_z)=\{1,\, 19,\, 43,\, 94\},\ E(2x^2+5y^2+T_z)=\{4,\, 27\},
\\&E(3x^2+4y^2+T_z)=\{2,\, 11,\, 23,\, 50,\, 116,\, 135,\, 138\},
\\&E(x^2+2y^2+6T_z)=\{5,\, 13,\, 46,\, 161\},
\\&E(x^2+5y^2+3T_z)=\{2,\, 11,\ 26,\, 37,\, 40,\,53,\,62,\, 142,\,220,\,425,\, 692\},
\\&\max E(2x^2+2y^2+5T_z)=2748,\ \max E(2x^2+4y^2+3T_z)=3185,
\\&\max E(x^2+4y^2+5T_z)=2352,\  \max E(2x^2+3y^2+5T_z)=933.
\end{align*}
Under the generalized Riemann hypothesis, the argument of Ono and
Soundararajan \cite{OnoSound} would allow one to use Waldspurger's
theorem \cite{Waldspurger} (or a Kohnen-Zagier variant \cite{Kohnen}
when the corresponding modular form is in Kohnen's plus space)
 to determine effectively a computationally feasible bound beyond which every integer
 is represented and hence verify that the above lists (and all lists contained herein) are indeed complete.
 This is done by carefully comparing the growth of the class numbers of imaginary quadratic fields with
 the growth of coefficients of a particular cusp form.

\medskip
Recall that  $\{x^2+2T_y:\ x,y\in\Z\}=\{T_x+T_y:\ x,y\in\Z\}$ as
observed by Euler. (See, e.g., (3.6.3) of \cite[p.71]{Berndt}, and
\cite[Lemma 1]{Sun1}.)
 Thus we say that $x^2+2T_y$ is equivalent to $T_x+T_y$ and denote this by
 $x^2+2T_y\sim T_x+T_y$.
 \medskip

 \begin{corollary}\label{abc}
Let $a,b,c\in\Z^+$ with $c$ odd. Then,
\begin{align*}&\text{all sufficiently large odd integers have the form}
\ 2ax^2+2by^2+cz^2
\\\iff&ax^2+by^2+4cT_z\ \text{is almost universal}
\\\iff&2\|ab,\ -ab\ R\ c,\ -2ac\ R\ b'\ \text{and}\ -2bc\ R\ a'.
\end{align*}
In particular,
\begin{align*}&\text{all
sufficiently large odd numbers are represented by}\
2ax^2+c(y^2+z^2)
\\\iff&ax^2+2cy^2+4cT_z\sim ax^2+2c(T_y+T_z)\ \text{is almost universal}
\\\iff& c=1,\ \text{and all prime divisors of }a\ \text{are congruent to 1 mod 4}.
\end{align*}
\end{corollary}

\medskip
\noindent{\it Remark} 1.3. In 2005 L. Panaitopol \cite{P} showed that
for $a,b,c\in\Z^+$ with $a\le b\le c$ all positive odd integers can be written
in the form $ax^2+by^2+cz^2$ with $x,y,z\in\Z$, if and only if the vector $(a,b,c)$
is $(1,1,2)$ or $(1,2,3)$ or $(1,2,4)$. For $n\in\N$, clearly $2n+1=x^2+2y^2+3z^2$ for some $x,y,z\in\Z$
if and only if there are $x,y,z\in\Z$ such that $2n+1=(8T_x+1)+2y^2+3(2z)^2$
(i.e., $n=4T_x+y^2+6z^2$) or
$2n+1=(2x)^2+2y^2+3(8T_z+1)$ (i.e., $n-1=2x^2+y^2+12T_z$). By Corollary \ref{abc},
the forms $6x^2+y^2+4T_z$ and $2x^2+y^2+12T_z$ are almost universal. Our computation suggests that
$$E(6x^2+y^2+4T_z)=\{2,\,3,\,17,\,23,\,38,\,51,\,86,\,188\}$$
and
$$E(2x^2+y^2+12T_z)=\{5,\,7,\,10,\,26,\,35,\,65,\,92,\,127,\,322\}.$$
\medskip

\begin{corollary}\label{ab}
Let $a,b\in\Z^+$ with $b$ odd. If $\SF(a')$ or $\SF(b)$ has a prime divisor $p\eq3\ (\mo\ 4)$
(which happens when $a'$ or $b$ is congruent to $3$ mod $4$),
then
\begin{align*}&ax^2+by^2+2T_z\ \text{is almost universal}
\\\iff&ax^2+y^2+2bT_z\ \text{is almost universal}
\\\iff&-a\R b\ \text{and}\ -b\ R\ a',
\end{align*}
and
\begin{align*}&ax^2+2y^2+bT_z\ \text{is almost universal}
\\\iff&ax^2+2by^2+T_z\ \text{is almost universal}
\\\iff& -2a\R\ b\ \text{and}\ -b\ R\ a'.
\end{align*}
\end{corollary}

\begin{corollary}\label{precisecor}
Let $a$ be any positive integer.

{\rm (i)} The form $ax^2+y^2+T_z$ is almost universal if and only if $-2\R a'$
(i.e., every odd prime divisor of $a$ is congruent to $1$ or $3$ modulo 8).
Also, $ax^2+2y^2+2T_z$ is almost universal if and only if
each prime divisor of $a$ is congruent to $1$ or $3$ modulo $8$.

{\rm (ii)}
\begin{align*}&ax^2+2y^2+T_z\ \text{is almost universal}
\\\iff &ax^2+y^2+2T_z\sim ax^2+T_y+T_z\ \text{is almost universal}
\\\iff&-1\R a',\ \text{i.e., every odd prime divisor of }a\ \text{is congruent to }1\ \mo\ 4.
\end{align*}
Also,
\begin{align*}& ax^2+2y^2+4T_z\sim ax^2+2T_y+2T_z\ \text{is almost universal}
\\\iff&\text{all prime divisors of }a\ \text{are congruent to }1\ \mo\ 4,
\end{align*}
and
\begin{align*}& ax^2+4y^2+2T_z\ \text{is almost universal}
\\\iff&a\eq 1\ (\mo\ 8)\ \text{and each prime divisor of }a\ \text{is congruent to }1\ \mo\ 4.
\end{align*}
\end{corollary}

\medskip
{\it Example} 1.2. By Corollary \ref{precisecor}, the form $5x^2+4y^2+2T_z$ is not almost universal
though it is asymptotically universal.
Also, our computation suggests the following information:
$$E(11x^2+y^2+T_z)=\{8,\,34,\,348\}\ \ \text{and}\ \ E(12x^2+y^2+T_z)=\{8,\,20,\,146,\,275\}.$$

\begin{corollary}\label{ExtendedCorollary}
Let $a\in\Z^+$. Then
\begin{align*}&ax^2+3y^2+T_z\ \text{is almost universal (or asymptotically universal)}
\\\iff&a\equiv 1\ (\mo\ 3),\ \text{and}\ \lfloor p/12\rfloor\ \text{is even for any odd prime divisor}
\ p\ \text{of}\ a,
\end{align*}
and
\begin{align*}&ax^2 +y^2 +3T_z\ \text{is almost universal (or asymptotically universal)}
\\\iff&a\equiv 2\ (\mo\ 3),\ \text{and}\ \lfloor p/12\rfloor\ \text{is even for any odd prime divisor}
\ p\ \text{of}\ a.
\end{align*}
Also,
\begin{align*}&ax^2+2y^2+6T_z\ \text{is almost universal (or asymptotically universal)}
\\\iff&a\equiv 1\ (\mo\ 6),\ \text{and}\ \lfloor p/12\rfloor\ \text{is even for any prime divisor}\ p\ \text{of}\ a,
\end{align*}
and
\begin{align*}&ax^2 +6y^2 +2T_z\ \text{is almost universal (or asymptotically universal)}
\\\iff&a\equiv 5\ (\mo\ 6),\ \text{and}\ \lfloor p/12\rfloor\ \text{is even for any prime divisor}\ p\ \text{of}\ a.
\end{align*}
\end{corollary}

\begin{corollary}\label{squarefree} Let $m$ be a positive integer.

{\rm (i)} $x^2+y^2+mT_z$ is almost universal if and only if $4\nmid
m$ and all odd prime divisors of $m$ are congruent to $1$ mod $4$.
Also, $2x^2+y^2+mT_z$ is almost universal if and only if $8\nmid m$
and each odd prime divisor of $m$ is congruent to $1$ or $3$ mod
$8$.

{\rm (ii)} Let $k\in\Z^+$. Then the form $2^{2k}x^2+y^2+mT_z$
is almost universal if and only if $4\nmid m$,  $-1\R m'$, and
$$2\| m\ \Longrightarrow\ m\ \text{is squarefree}.$$
Also, the form $2^{2k+1}x^2+y^2+mT_z$ is almost universal if and only if
$4\nmid m$, $-2\R m'$, and
$$m\eq1\ (\mo\ 8)\ \ \Longrightarrow\ \ m\ \text{is squarefree}.$$

{\rm (iii)} Let $k,l\in\Z^+$ with $k\ge l$. Then $2^kx^2+2^ly^2+mT_z$ is asymptotically universal
if and only if for each prime divisor $p$ of $m$ we have
$$\begin{cases} p\eq1\ (\mo\ 4)&\text{if}\ k\eq l\ (\mo\ 2),
\\p\eq1\ \text{or}\ 3\ (\mo\ 8)&\text{otherwise}.
\end{cases}$$
When $2^kx^2+2^ly^2+mT_z$ is asymptotically universal, it is almost universal if and only if
$m$ is squarefree, or both $2\mid k$ and $l=1$.
\end{corollary}

\medskip
{\it Example} 1.3. By Corollary \ref{squarefree}, the forms $4x^2+y^2+50T_z$, $8x^2+y^2+9T_z$
and $2x^2+2y^2+25T_z$ are not almost universal though they are asymptotically universal. We also have
the following observations via computation:
\begin{align*}&\max E(x^2+y^2+10T_z)=546,\ \max E(2x^2+y^2+11T_z)=985;
\\&\max E(4x^2+y^2+10T_z)=5496,\ \max E(4x^2+2y^2+9T_z)=9555;
\\&\max E(2x^2+2y^2+13T_z)=22176,\ \max E(8x^2+y^2+3T_z)=499.
\end{align*}

\begin{corollary}\label{ExtendedCorollary2} Let $a\in\integer^{+}$ be even.

{\rm (i)} Suppose that $v_2(a)$ is even and each odd prime divisor
of $a$ is congruent to $1$ modulo $3$.  Then $ax^2+216 y^2 + T_z$ is
asymptotically universal. Moreover, this form is not almost
universal if and only if every prime divisor of $\SF(a')$ is
congruent to $1$ or $19$ modulo $24$, and the number of prime
divisors congruent to $19$ modulo $24$ is odd.

{\rm (ii)} Assume that $v_2(a)$ is odd, $a'\eq\pm1\ (\mo\ 10)$ and
$2\mid\lfloor p/10\rfloor$ for every prime divisor $p$ of $a'$.
Then $ax^2+250y^2 + T_z$ is asymptotically universal. Moreover,
this form is not almost universal if and only if $a'\eq 21,29\
(\mo\ 40)$ and every prime divisor of $\SF(a')$ is congruent to
$1$ or $9$ modulo $20$.
\end{corollary}

\medskip
\noindent {\it Remark} 1.4. Corollary \ref{ExtendedCorollary2}
implies that the forms $76x^2+216y^2+T_z$ and $58x^2+250y^2+T_z$
are asymptotically universal but not almost universal.
\medskip

For the form $ax^2+bT_y+cT_z$, we obtain the following result.

\begin{theorem}\label{onesufflargenec}
Let $a,b,c\in \integer^{+}$ with $\gcd(a,b,c)=1$ and $v_2(b)\ge v_2(c)$.
Consider the form
$$f(x,y,z):=ax^2+bT_y + cT_z
$$
and assume that both $(1)$ and $(2)$ in Theorem \ref{onealmost} hold.

{\rm (i)} When $v_2(b)\not\in\{3,4\}$,
$f$ is not almost universal if and only if we have the following $(1)-(4)$.

\begin{enumerate}
 \item\label{onecongcond} $4\nmid b+c$ and $\SF(a'b'c')\equiv (b+c)' \pmod{2^{3-v}}$, where $v:=v_2(b+c)<2$.
 \item \label{oneprimeprodcond} All prime divisors of $\SF(a'b'c')$ are congruent to $1$ or $3$ modulo $8$
 if $\SF(abc)\equiv b+c\pmod{2}$, and congruent to $1$ modulo $4$ otherwise.
 \item \label{onerepcond} $8ax^2+by^2+cz^2 = 2^{v} \SF(a'b'c')$ has no integral solutions with $y$ and $z$ odd.
\item\label{oneiff}
$$\begin{cases}v_2(b)\le1\Rightarrow v_2(a)-v_2(b)\in\{2,4,6,\ldots\},&
\\v_2(b)=2\Rightarrow v_2(a)\in\{1,3,5,\ldots\},&
\\ v_2(b)\in\{5,7,\ldots\}\Rightarrow (4\mid a\ \text{or}\ 2\mid c),&
\\v_2(b)\in\{6,8,\ldots\}\Rightarrow (2\mid a\ \text{or}\ a\eq c\ (\mo\
8)).
\end{cases}$$
\end{enumerate}

{\rm (ii)} In the case $v_2(b)\in\{3,4\}$,  if $f$ is not almost universal, then
$(1)-(3)$ above hold and also
$$\begin{cases}v_2(b)=3\Rightarrow (4\mid a\ \text{or}\ 2\mid c),&
\\v_2(b)=4\Rightarrow (2\mid a\ \text{or}\ a\eq c\ (\mo\ 8)).
\end{cases}$$
Moreover, provided $(1)-(3)$ in part {\rm (i)} and the condition $2\nmid v_2(a)$,
$f$ is not almost universal if $v_2(b)=4$, or $v_2(a)\ge v_2(b)=3$ and $b'\eq c'\ (\mo\ 8)$.
\end{theorem}

\medskip

{\it Example 1.4}. Consider those forms $ax^2+bT_y+cT_z$ with
$a,b,c\in\Z^+$ and $a+b+c\le10$. By Theorem \ref{onesufflargenec},
we find that those asymptotically universal ones are almost
universal. Below is a complete list of those forms $ax^2+bT_y+cT_z$
with $a,b,c\in\Z^+$ and $a+b+c\le10$ which are almost universal but
not universal:
\begin{align*}&5x^2+T_y+T_z\sim x^2+5y^2+2T_z,\ 5x^2+2T_y+2T_z\sim 2x^2+5y^2+4T_z,
\\&8x^2+T_y+T_z\sim x^2+8y^2+2T_z,
\ 2x^2+3T_y+2T_z,\ x^2+4T_y+3T_z,
\\&2x^2+5T_y+T_z,\  4x^2+3T_y+T_z,
\ 3x^2+5T_y+T_z,\ 3x^2+4T_y+2T_z,
\\&4x^2+4T_y+T_z,\ 6x^2+2T_y+T_z,\ 5x^2+3T_y+2T_z,\ 5x^2+4T_y+T_z.
\end{align*}
For the above forms from the second line, our computation via
computer suggests the following information:
\begin{align*}&E(8x^2+T_y+T_z)=E(x^2+8y^2+2T_z)=\{5,\, 40,\, 217\},
\\&E(2x^2+3T_y+2T_z)=\{1,\,16\},\ E(x^2+4T_y+3T_z)=\{2,\,6,\,80\},
\\&E(2x^2+5T_y+T_z)=\{4\},\ E(4x^2+3T_y+T_z)=\{2,\,11,\,27,\,38,\,86,\,93,\,188,\,323\},
\\&E(3x^2+5T_y+T_z)=\{2,\,7\},\ E(3x^2+4T_y+2T_z)=\{1,\,8,\,11,\,25\},
\\&E(4x^2+4T_y+T_z)=\{2,\,108\},\ E(6x^2+2T_y+T_z)=\{4\},
\\&E(5x^2+3T_y+2T_z)=\{1,4,13,19,27,46,73,97,111,123,151,168\},
\\&E(5x^2+4T_y+T_z)=\{2,\,16,\,31\}.
\end{align*}

In Corollary \ref{precisecor} we determined when $ax^2+T_x+T_y$ or $ax^2+2T_x+2T_y$
is almost universal. The following corollary deals with two other similar forms.

\begin{corollary}\label{Tx+2Ty} Let $a$ be a positive integer. Then $ax^2+2T_y+T_z$ is almost universal
if and only if all odd prime divisors of $a$ are congruent to $1$ or $3$ mod $8$. Also,
$ax^2+4T_y+T_z$ is almost universal if and only if all odd prime divisors of $a$ are congruent to $1$ mod $4$.
\end{corollary}

\medskip
{\it Example} 1.5. By means of computation, we believe that
$$E(9x^2+2T_y+T_z)=\{4\}\ \ \text{and}\ \ E(11x^2+2T_y+T_z)=\{4,\,25,\,94,\,123\}.$$
\medskip

\begin{corollary}\label{x^2+Ty} Let $m$ be any positive integer.

{\rm (i)} If all odd prime divisors of $m$ are congruent to $1$ or $3$ mod $8$, and
$m'\eq 3\ (\mo\ 8)$ or $v_2(m)\not=4,6,\ldots$, then
$x^2+T_y+mT_z$ is almost universal. The converse also holds when $v_2(m)\not=4$.

{\rm (ii)} For $k\in\Z^+\setminus\{3,4\}$, the form $2^k(x^2+T_y)+mT_z$
is almost universal if and only if $k\in\{1,2\}$
and all prime divisors of $m$ are congruent to $1$ or $3$ mod $8$.
When $m\eq1\ (\mo\ 8)$, the form $8(x^2+T_y)+mT_z$ is not almost universal.
\end{corollary}

{\it Example}\ 1.6.
By Corollary \ref{x^2+Ty}, the form $8x^2+8T_y+T_z$ is not almost universal though it is asymptotically universal.
We also have the following guess based on our computation:
$$E(x^2+T_y+9T_z)=\{8,\,47\},\ E(x^2+T_y+11T_z)=\{8\},
\ E(x^2+T_y+12T_z)=\{8,\,20\}.$$

\medskip

\begin{corollary}\label{Tx+Ty} Let $m$ be any positive integer.

{\rm (i)} When $v_2(m)\not=3$, the form $x^2+2T_y+mT_z\sim T_x+T_y+mT_z$
is almost universal if and only if all odd prime divisors of $m$
are congruent to $1$ mod $4$ and $v_2(m)\not=5,7,\ldots$.

{\rm (ii)} For $k\in\Z^+\setminus\{2\}$, the form
$2^k(x^2+2T_y)+mT_z\sim 2^k(T_x+T_y)+mT_z$ is almost universal if and only if
$k=1$ and all prime divisors of $m$ are congruent to $1$ mod $4$.
\end{corollary}

\noindent{\it Remark}\ 1.5.
By Corollary \ref{Tx+Ty}(ii), $8x^2+16T_y+T_z\sim 8(T_x+T_y)+T_z$ is not almost universal
though it is asymptotically universal.
In \cite{Sun4} the second author conjectured that any integer $n>1029$ is either a triangular number
or a sum of two odd squares and a triangular number (i.e., $n=(8T_x+1)+(8T_y+1)+T_z$ for some $x,y,z\in\Z$);
in other words,
$$E(8T_x+8T_y+T_z)\cap[1028,+\infty)\subseteq\{T_m-2:\ m\in\Z^+\}.$$
Recently, Oh and the second author \cite{Sun2} showed that $T_m-2\in E(8T_x+8T_y+T_z)$
(i.e., $T_m$ is not a sum of two odd squares and a triangular number) if and only if
$2m+1$ is a prime congruent to 3 mod 4.
\medskip

{\it Example} 1.7. Via computation we make the following observation:
\begin{align*}&E(x^2+2T_y+10T_z)=E(T_x+T_y+10T_z)=\{5,\,8\},
\\&E(x^2+2T_y+13T_z)=E(T_x+T_y+13T_z)=\{5,\,8,\,32,\,53\}.
\end{align*}

\begin{theorem}\label{trisufflarge}
Let $a,b,c\in\Z^+$ with $v_2(a)\ge v_2(b)\ge v_2(c)=0$. Assume that $-bc\R a'$, $-ac\R b'$ and $-ab \R c'$.
Consider  the form
$$
f(x,y,z):=aT_x+bT_y + cT_z.
$$

{\rm (i)} If $f$ is not almost universal, then we have the following $(1)-(4)$.

\begin{enumerate}
 \item\label{tricongcond} $4\nmid a+b+c$ and $\SF(a'b'c')\equiv (a+b+c)' \pmod{2^{3-v}}$, where $v=v_2(a+b+c)<2$.
 \item \label{triprimeprodcond}
 All prime divisors of $\SF(a'b'c')$ are congruent to $1$ modulo $4$ if $\SF(abc)\equiv a+b+c\ (\mo\ 2)$,
 and congruent to $1$ or $3$ modulo $8$ otherwise.
 \item \label{trirepcond} $ax^2+by^2+cz^2 = 2^{v} \SF(a'b'c')$ has no integral solution with $x,y,z$ all odd.
\item \label{trinec}
$$\begin{cases}v_2(b)\leq 1\Rightarrow v_2(a)-v_2(b)\in\{3,5,7,\ldots\},&
\\v_2(b)=2\Rightarrow v_2(a)\in\{2,4,6,\ldots\}.
\end{cases}$$
\end{enumerate}

{\rm (ii)} $f$ is not almost universal under  $(\ref{tricongcond})-(\ref{trirepcond})$ in part {\rm (i)},
and the following condition stronger than $(\ref {trinec})$:
$$\begin{cases}v_2(b)\leq 1\Rightarrow v_2(a)-v_2(b)\in\{5,7,\ldots\}&
\\v_2(b)\in\{2,4\}\Rightarrow v_2(a)\in\{4,6,\ldots\}&
\\v_2(b)=3\Rightarrow (v_2(a)\in\{6,8,\ldots\}\ \&\ b'\eq c\ (\mo\ 8)).
\end{cases}$$
\end{theorem}

\medskip

{\it Example} 1.8. Consider those forms $aT_x+bT_y+cT_z$ with
$a,b,c\in\Z^+$ and $a+b+c\le10$. By Theorem \ref{trialmost}, we find
the following complete list of those asymptotically universal forms
which are not universal:
\begin{align*}&T_x+4T_y+4T_z\sim 4x^2+8T_y+T_z,
\\&T_x+T_y+8T_z\sim x^2+8T_y+2T_z,
\\&2T_x+2T_y+5T_z\sim 2x^2+4T_y+5T_z,
\\&T_x+2T_y+6T_z,\ 2T_x+3T_y+4T_z,\  T_x+4T_y+5T_z.
\end{align*}
By Theorem \ref{trisufflarge}, the last four forms are in fact almost universal;
our computation via computer
suggests the following information:
$$E(2T_x+2T_y+5T_z)=E(2x^2+4T_y+5T_z)
=\{1,3,10,16,28,43,46,85,169,175,211,223\}$$ and
$$E(T_x+2T_y+6T_z)=\{4,\,50\},\ E(2T_x+3T_y+4T_z)=\{1,\, 8,\,
31\}, \ E(T_x+4T_y+5T_z)=\{2\}.$$ As for the first two forms
$T_x+4T_y+4T_z$ and $T_x+T_y+8T_z$, neither Theorem
\ref{trisufflarge} nor Theorem \ref{onesufflargenec} tells us
whether they are almost universal or not. However, with some special
arguments, the first author \cite{Kane4} was able to show that they
are not almost universal though they are asymptotically universal.
By Theorem 1.1(ii) of an earlier paper \cite{Sun2},
$E(T_x+T_y+8T_z)$ actually consists of those $2T_m-1$ ($m\in\Z^+$)
with $2m+1$ having no prime divisors congruent to $3$ mod $4$;
similarly, by \cite[Theorem 1(iii)]{Sun1} and \cite[Theorem
2.1(ii)]{Sun2}, $E(T_x+4T_y+4T_z)$ consists of those $T_m-1$
($m\in\Z^+$) with $2m+1$ having no prime divisors congruent to $3$
mod $4$.
\medskip

\begin{corollary}\label{TTT} Let $a\in\Z^+$. Then the form $aT_x+2T_y+T_z$ is
almost universal if  each odd prime divisor of $a$ is congruent to
$1$ or $3$ mod $8$, and either $a'\eq1\ (\mo\ 8)$ or
$v_2(a)\not=4,6,\ldots$. We also have the converse when
$v_2(a)\not=4$.
\end{corollary}

\medskip
\noindent{\it Remark} 1.6. In \cite{Kane4} the first author was
able to show that the special form $48T_x+2T_y+T_z$ is not almost
universal (though it is asymptotically universal by Theorem
\ref{trialmost}).

\medskip
{\it Example} 1.9.  Our computation leads us to make the following
observation:
\begin{gather*} E(9T_x+2T_y+T_z)=\{4,\,46\},\ E(11T_x+2T_y+T_z)=\{4,\,25\},
\\E(22T_x+2T_y+T_z)=\{4,\,11,\,14,\,19,\,46,\,54\}.
\end{gather*}
\smallskip

Our following conjecture is a supplement to Theorems \ref{onesufflargenec}
and \ref{trisufflarge}; its solution might involve a
 further investigation of the spinor norm mapping or alternation of certain coefficients of cusp
 forms.

\begin{conjecture}\label{remaining} Let $a,b,c$ be positive
integers.

{\rm (i)} In the case $v_2(b)\ge v_2(c)$ and $v_2(b)\in\{3,4\}$, if  $(1)-(3)$ in Theorem
$\ref{onesufflargenec}$ hold and also
$$\begin{cases}v_2(b)=3\Rightarrow (4\mid a\ \text{or}\ 2\mid c)&
\\v_2(b)=4\Rightarrow (2\mid a\ \text{or}\ a\eq c\ (\mo\ 8)),
\end{cases}$$
then the form $ax^2+bT_y+cT_z$ is not almost universal.

{\rm (ii)} In the case $v_2(a)\ge v_2(b)\ge v_2(c)=0$, if $(1)-(3)$ in Theorem
$\ref{trisufflarge}$ hold, and
$$v_2(a)=v_2(b)=2\ \ \text{or}\ \ v_2(a)=v_2(b)+3\in\{3,4\}\ \ \text{or}\ \ v_2(b)\in\{3,4\},$$
then the form $aT_x+bT_y+cT_z$ is not almost universal.
\end{conjecture}

In the next section we are going to introduce some further notation
and give an overview of our method. In Section 3 we will deal with
asymptotically universal forms and prove Theorems 1.1-1.3 and
Corollaries 1.4-1.5. Section 4 is devoted to the proofs of the
remaining theorems and corollaries concerning almost universal
forms.

\section{Notation and Brief Overview}\label{NotationSection}
Our arguments will involve the theory of modular forms and spinor
exceptional square classes for quadratic forms. A good introduction
to modular forms may be found in Ono's book \cite{Ono1} and a good
introduction to quadratic forms may be found in O'Meara's book
\cite{Omeara1}.  We will first reduce the questions at hand to
questions about certain related (ternary) quadratic forms.  Since
$8T_x+1$ is an odd square, multiplying by $8$ and adding some
positive integer will give a form $Q(x,y,z)$ which is a sum of
squares  with the restriction that some of $x$, $y$, and $z$ must be
odd. If we take $r_Q(n)$ to be the number of solutions to
$Q(x,y,z)=n$ with the given restrictions on $x,y,z\in\integer$, then
define
$$
\theta_Q(\tau):=\sum_{n=0}^{\infty} r_Q(n) q^n,
$$
where $q=e^{2\pi i\tau}$ with $\tau$ in the upper half plane.  Since
the number of solutions with $z$ odd equals the number of solutions
with $z$ arbitrary minus the number of solutions with $z$ even and
since the form with $z$ even gives another quadratic form, we get an
inclusion/exclusion of theta series of
 quadratic forms.  Let a ternary quadratic form $Q'(x,y,z)$ be given.  Then it is well known that the theta series
$$
\theta_{Q'}(\tau):= \sum_{n=0}^{\infty} r_{Q'}(n) q^n
$$
is a modular form of weight $3/2$, where $r_{Q'}(n)$ is the number
of solutions to $Q'(x,y,z)=n$ with $x,y,z\in \integer$. The theta
series splits naturally into the following three parts
$$
\theta_{Q'}=\theta_{gen(Q')} +
\left(\theta_{spn(Q')}-\theta_{gen(Q')}\right) +
\left(\theta_{Q'}-\theta_{spn(Q')}\right),
$$
where the $n$-th coefficients of $\theta_{gen(Q')}$ and
$\theta_{spn(Q')}$ are the weighted average of the number of
representations of $n$ by the genus and the spinor genus of $Q'$,
respectively. Furthermore, $\theta_{gen(Q')}$ is an Eisenstein
series, $\theta_{spn(Q')}-\theta_{gen(Q')}$ is a cusp form in the
space of one dimensional theta series, and
$\theta_{Q'}-\theta_{spn(Q')}$ is a cusp form in the orthogonal
complement of the space of one dimensional theta series.  For a full
description, see the survey paper of Schulze-Pillot
\cite{SchulzePillot1}. We will then use the argument of Duke and
Schulze-Pillot \cite{DSP}.

The coefficients of $\theta_{spn(Q')}-\theta_{gen(Q')}$ are supported at finitely many square classes.
If $r_{Q',p^k}(n)$ is the number of solutions to $Q'(x,y,z)=n\pmod{p^k}$ with $x,y,z\in \integer/p^k\integer$,
then the $n$-th coefficient of the Eisenstein series was shown by Siegel (cf. Jones \cite{Jones}) to be
\begin{equation}\label{EisGrowthEqn}
\prod_{p\text{ prime}}\lim_{k\to\infty} \frac{r_{Q',p^k}(n)}{p^{2k}}.
\end{equation}
An anisotropic prime $p$ is a prime for which the equation $Q'(x,y,z)=0$ has no non-trivial solutions in $\integer_p$.
Notice that for $q\neq p$, $r_{Q',q^k}(np^2)=r_{Q',q^k}(n)$, since $p$ is invertible and hence we have a bijection
between solutions to $Q'(x,y,z)=np^2$ and $Q'(x',y',z')=n$ by taking $(x,y,z)\to p^{-1}(x,y,z)$.
Because $Q'(x,y,z)=0$, if $n$ has sufficiently large divisibility by $p$ (i.e., the $p$-adic order of $n$ is sufficiently large), then it is easy to check that
$r_{Q',p^k}(np^2)=r_{Q',p^k}(n)$, and hence the $np^{2k}$ coefficients of the Eisenstein series grow like a constant
with respect to $k$.

When $n$ has bounded divisibility at every anisotropic prime (i.e., the orders of $n$ at anisotropic primes are bounded), Equation (\ref{EisGrowthEqn}),
and hence the coefficients of the Eisenstein series, grow like a certain class number
(cf. Jones \cite[Theorem 86]{Jones}), and hence are (ineffectively) $\gg n^{1/2-\epsilon}$
by the bound of Siegel \cite{Siegel}.  The coefficients of the cusp forms in the orthogonal complement of
one dimensional theta series (to which $\theta_{Q'}-\theta_{spn(Q')}$ belongs) were
$\ll n^{1/2-1/28+\epsilon}$ as first obtained by Duke \cite{Duke},
extending Iwaniec's result \cite{Iwaniec} (for coefficients of half integral weight $\ge 5/2$ modular forms)
to the case of weight $3/2$ modular forms.  Better bounds have been obtained by the amplification method
on sums of special values of $L$-series as in Blomer, Harcos, and Michel \cite{BHM}, but the bound above is sufficient
to guarantee that if $n$ has bounded divisibility at the anisotropic primes and $n$ is not in one of the finitely many
square classes where the coefficients of $\theta_{spn(Q')}-\theta_{gen(Q')}$ are supported, the growth of
the coefficients of the Eisenstein series $\theta_{gen(Q')}$ will overwhelm the growth of the coefficients of
the cusp form $\theta_{Q'}-\theta_{gen(Q')}$ and hence the coefficients of $\theta_{Q'}$ will be positive for
sufficiently large $n$, with bounded divisibility by the anisotropic primes, outside of these finitely many square
classes.  Moreover, if we take a weighted sum
$$
\sum_{i=1}^m w_i \theta_{Q_i'},
$$
such as the inclusion/exclusion above, of finitely many such $\theta_{Q_i'}$ where
$\theta_{gen(Q_i')}=c_i\theta_{gen(Q')}$ and $\sum_{i=1}^{m} w_i c_i>0$, then the resulting theta series will be
$$\bigg(\sum_{i=1}^m w_i c_i\bigg) \theta_{gen(Q')}+f_1+f_2,
$$
where $f_1=\sum_{i=1}^m w_i(\theta_{spn(Q_i')}-\theta_{gen(Q_i')})$ and
$f_2=\sum_{i=1}^{m} w_i (\theta_{Q_i'}- \theta_{spn(Q_i')})$.
The bound of Duke \cite{Duke} given above then shows that outside of the finitely many square classes
where the coefficients of $f_1$ are supported the $n$-th coefficient of this weighted average is positive
for sufficiently large $n$ with bounded divisibility at the anisotropic primes.
The condition of the bounded divisibility at anisotropic primes will pose only a minor complication, and
we will find in the end that for any asymptotically universal form the associated quadratic forms will never have an
anisotropic prime $p\neq 2$, while conditions modulo 8 will guarantee that the coefficients which we are interested
in automatically have bounded divisibility by $2$.

We will thus be interested in determining which square classes of
coefficients $t\integer^2$ are supported by
$\theta_{spn(Q')}-\theta_{gen(Q')}$.  Kneser \cite{Kneser1} gave a
necessary condition and later Schulze-Pillot \cite{SchulzePillot2}
extended this to give necessary and sufficient conditions. For a
quadratic form $Q'$, there is an associated bilinear form
$B(x,y)=\left(Q'(x+y)-Q'(x)-Q'(y)\right)/2$. We will call $V$ a
(ternary) quadratic space over $\rational_2$ if it is a finite
dimensional vector space over $\rational_2$ with an associated
bilinear form $B$.  There is a quadratic form (over $\rational_2$)
associated to $V$ given by $Q'(x)=B(x,x)$ for every $x\in V$.  Fix a
$\integer_2$-lattice $L$.  The quadratic form (over $\integer_2$)
associated to $L$ is $Q'(x)=B(x,x)$ with $x\in L$. In our case, the
lattice will always have an orthogonal basis $x_1,x_2,x_3$ with
$B(x_i,x_j)=0$ when $i\neq j$.  We will denote the $\integer_2$
lattice whose corresponding quadratic form is $a x^2+b y^2 +cz^2$ by
$\left\langle a,b,c\right\rangle_2$.

We will denote isometries from $V$ to $V$ by $O(V)$. Let $O^+(V)$ be
the subgroup of rotations consisting of isometries with determinant
1.  We also use $O^+(L)$ to denote the rotations which fix $L$. Each
rotation is the product of an even number of symmetries, where the
symmetry $\tau_v$ with $v\in V$ is defined by
$$x\mapsto x-\frac{2B(x,v)}{Q'(v)} v.$$
The spinor norm mapping is the mapping
$\theta(\sigma)=Q'(v_1)\cdots Q'(v_m){\rational_2\cross}^2$ where
$\sigma =\tau_{v_1}\cdots \tau_{v_m}$.  The set $\theta(O^{+}(L))$
forms a subgroup of $\rational_2\cross/{\rational_2\cross}^2$. For
the $2$-adic lattice $L=L_2=\langle a,b,c\rangle_2$, Earnest and
Hsia determined this subgroup explicitly in \cite{EarnestHsia2}.

Fix an imaginary quadratic field extension $K/\rational$ (in our cases, $K=\rational(i)$ or $K=\rational(\sqrt{-2})$)
 and a prime ideal $\beta$ (of the ring $O_K$ of algebraic integers in $K$)
 dividing $2$.  For convenience, we define
$$
K_n:=\rational(\sqrt{-n\SF(n)}).
$$
We say that $\alpha\in \rational_2\cross /{\rational_2\cross}^2$ is
a local norm at $2$ (from the completion $K_{\beta}$ to
$\rational_2$) if $\alpha=x^2+ny^2$ for some $x,y\in \rational_2$.
We will denote the set of local norms at $2$ by $N_2(K)$.
 Note that $$N_2(\rational(i)) = {\rational_2\cross}^2\cup 5{\rational_2\cross}^2\cup 2{\rational_2\cross}^2
 \cup 10{\rational_2\cross}^2$$ and
  $$N_2(\rational(\sqrt{-2}))={\rational_2\cross}^2\cup 3{\rational_2\cross}^2
 \cup 2{\rational_2\cross}^2\cup 6{\rational_2\cross}^2.$$

Using explicit results of Earnest, Hsia, and Hung
\cite{EarnestHsia1} based on Schulze-Pillot's classification of
spinor exceptional square classes \cite{SchulzePillot2}, we will
reduce the question at hand to showing Kneser's necessary condition
at the prime $2$. The necessary condition of Kneser which we will
need to show is that $\theta(O^+(L))\subseteq N_2(K)$ (cf.
\cite{Kneser1}). We will use the explicit results of Earnest and
Hsia \cite{EarnestHsia2} to determine when the necessary condition
is satisfied.

 For $a,b\in\rational_2\cross $, the Hilbert symbol
$(a,b)_2\in\{\pm1\}$ takes the value 1 if and only if
$ax^2+by^2=z^2$ for some $x,y,z\in\rational_2$ with $x,y,z$ not all zero.
We will need the following theorems.

\begin{theorem}[Earnest and Hsia \cite{EarnestHsia2}]\label{EarnestHsiaTheorem}
 Let $U$ denote the group of units in $\integer_2$ and let $\alpha\in U$. Then
$$
\theta(O^{+}(\left\langle1,2^r \alpha\right\rangle_2)) =
\left\{
\begin{array}{ll}
\{ \gamma\in \rational_2\cross : (\gamma, -2\alpha)_2=1\}& \text{if } r\in\{1,3\},\\
\{ \gamma\in U{\rational_2\cross}^2 : (\gamma, -\alpha)_2=1\}& \text{if } r=2,\\
{\rational_2\cross}^2 \cup \alpha {\rational_2\cross}^2 \cup 5{\rational_2\cross}^2 \cup 5\alpha {\rational_2\cross}^2
& \text{if }r=4,\\
{\rational_2\cross}^2 \cup \alpha {\rational_2\cross}^2& \text{if }r\geq 5.
\end{array}
\right.
$$
\end{theorem}

Furthermore, Earnest and Hsia \cite[Theorem 2.2]{EarnestHsia2} showed  that for the lattice $L_2:=\langle c',2^rb',
2^sa'\rangle_2$, we have $\theta(O^+(L_2))=\rational_2\cross$ if $\{r,s-r\}\cap\{1,3\}\not=\emptyset$
and $\{r,s,s-r\}\cap\{2,4\}\not=\emptyset$.  If
$0<r<s$ and the conditions of Theorem 2.2 in \cite{EarnestHsia2}
are not satisfied, they proved that $\theta(O^+(L_2))$ is equal to
the union of the spinor norm restricted to 2-dimensional
sublattices, allowing us to use the above theorem.  Moreover, if
$s\geq 5$ and $r\in\{0,s\}$ then their argument follows
mutatis mutandis and will also allow us to reduce the problem to 2-dimensional sublattices.

Since our base field is $\rational_2$ and $K_{\beta}/\rational_2$ is ramified for $K=\rational(i)$
and $K=\rational(\sqrt{-2})$ we will only need the following restriction of the $2$-adic conditions from Earnest, Hsia,
and Hung's theorem.
\begin{theorem}[Earnest, Hsia, and Hung \cite{EarnestHsia1}]\label{EarnestHsiaHungTheorem}
Let $a,b,c\in\Z^+$ and $K=\rational(\sqrt{-abc})$.
Let $L_2=\left\langle c',2^rb',2^sa'\right\rangle_2$ with $0\leq r\leq s$, and let $t\in\Z^+$.
Assume that $\theta(O^+(L_2))\subseteq N_2(K)$, and define
$$L'=\begin{cases}\left\langle 2^{r-2}c', 2^r b', 2^s a'\right\rangle_2&\text{if}\ r+s\eq v_2(t)\ (\mo\ 2),
\\\left\langle 2^{r-3}c', 2^r b', 2^s a'\right\rangle_2&\text{otherwise}.\end{cases}$$
Consider the necessary and sufficient conditions given by Schulze-Pillot \cite{SchulzePillot2} for
 $t$ to be a primitive spinor exception for the genus of the quadratic form $Q(x,y,z)=ax^2+by^2+cz^2$.

\begin{enumerate}
\item  Set $L'':=\left\langle2^r c', 2^r b', 2^s a'\right\rangle_2$. When $r+s\eq v_2(t)\ (\mo\ 2)$,
 the Schulze-Pillot conditions are not satisfied
if and only if one of the following $(a)-(d)$ holds.
    \begin{enumerate}
        \item $r$ is odd and $v_2(t)\geq r-3$.
        \item $r$ is even, $\theta(O^+(L'))\not\subseteq N_2(K)$, and
        $$(r\neq s\ \&\ v_2(t)\geq r-2)\ \ \text{or}\ \ (r=s\ \&\ v_2(t) \geq r).$$
        \item $r$ is even, $\theta(O^+(L'))\subseteq N_2(K)$,  $\theta(O^+(L''))\not \subseteq N_2(K)$ and $v_2(t)\geq r$.
        \item $r$ is even, $\theta(O^+(L'))\subseteq N_2(K)$,  $\theta(O^+(L''))\subseteq N_2(K)$ and $v_2(t)\geq s$.
    \end{enumerate}
\item  When $r+s\not\eq v_2(t)\ (\mo\ 2)$, we have $0<r<s$, and
the Schulze-Pillot conditions are not satisfied if and only if one of
the following $(a)-(c)$ holds:
    \begin{enumerate}
        \item $r$ is even and $v_2(t)\geq r-4$.
        \item $r$ is odd, $\theta(O^+(L'))\not\subseteq N_2(K)$ and $v_2(t)\geq r-3$.
        \item $r$ is odd, $\theta(O^+(L'))\subseteq N_2(K)$ and $v_2(t)\geq s-2$.
    \end{enumerate}
\end{enumerate}

\end{theorem}

\section{On Asymptotically Universal Forms}\label{AlmostSection}

In this section, we will show which forms are asymptotically universal, proving Theorems \ref{twoalmost},
\ref{onealmost}, and \ref{trialmost}.  We will first need the following lemma.

\begin{lemma}\label{plocallemma}
Fix $a,b,c\in\integer^{+}$ with $\gcd(a,b,c)=1$.  Then
$$
Q(x,y,z)= a x^2+b y^2+c z^2
$$
represents every integer $p$-adically for each odd prime $p$ if and
only if  we have
\begin{equation}\label{plocaleqn} -ab\R c', \ -ac\R b',\ \text{and}\ -bc \R a'.
\end{equation}
\end{lemma}
\begin{proof}
It is well known that $Q$ represents every integer $p$-adically for any odd prime $p$ not dividing $abc$.

Let $p$ be an odd prime divisor of $abc$. Without loss of generality, we assume that $p\mid c$.
If $p\mid ab$, then $Q$ clearly only represents all squares or all non-squares modulo $p$.
Therefore, $p\nmid ab$.

Let a unit $u\in \integer_p$ be given.  Suppose that there are $x,y\in\Z_p$ such that
$ax^2+b y^2 =pu$.
If $x\in p\Z_p$ , then we have $y\in p\Z_p$ and hence $u\in p\Z_p$, which contradicts the fact
that $u$ is a unit.
Therefore both $x$ and $y$ must be units.  Taking the Legendre symbol of both sides yields
that
$$
\left(\frac{b}{p}\right)=\left(\frac{b y^2}{p}\right)=\left(\frac{pu -ax^2}{p}\right)
=\left(\frac{-a x^2}{p}\right)
=\left(\frac{-a}{p}\right).
$$

Now assume that $(\frac{-ab}{p})\neq 1$.  Let a unit
$u\in\integer_p$ be given. From the above, we know that $ax^2+by^2
= up$ does not have a solution. Suppose that there are $x,y,z\in
\Z_p$ such that $ax^2+by^2+cz^2 = u p$.  Then
$$
ax^2+by^2=\bigg(u-\frac{c}{p}z^2\bigg)p.
$$
If $p^2\mid c$, then $u-\frac{c}{p}z^2$ is also a unit, and it
follows that $up$ is not represented.  In the case $p\| c$, without
loss of generality we assume that the unit $\frac{c}{p}$ is a
square. If $u$ is not a square, then $u-\frac{c}{p}z^2$ must also be
a unit and it follows that $up$ is not represented.

By the above, (\ref{plocaleqn}) is a necessary condition.

If $n\in \integer_p$ is represented then so is $np^2$. Thus we only
need to show that $(\f{-ab}p)=1$ implies that those $n\in\Z_p$ with
$v_p(n)\in\{0,1\}$ are $p$-adically represented. We have already
shown that $(\f{-ab}p)=1$ if and only if those $n\in\Z_p$ with
$v_p(n)=1$ are represented, so we only need to prove that every unit
is represented. Hence, it suffices to show that at least one square
and one non-square are represented by $Q$. We will prove that $ax^2+b y^2$
represents every integer $p$-adically.  If $-1$ is not a
square, then
$$
\left(\frac{b}{p}\right)=-\left(\frac{a }{p}\right),
$$
and hence both squares and non-squares are represented. So we may assume that $-1$ is a square.
For any unit $u\in\Z_p$, the form $ax^2+b y^2$ represents every integer $p$-adically
if and only if $uax^2+ ub y^2$ represents every integer
$p$-adically. So, without loss of generality we may suppose that
$$
\left(\frac{a}{p}\right)=\left(\frac{b }{p}\right)=1.
$$
Since $-1$ is a square and we represent all squares by $ax^2$ (and also $b y^2$), we must represent $-1$.
We now argue inductively by noting that if $-m$ is a square, then $-m-1$ is represented by $a x^2+b y^2$ via
taking $a x^2=-m$ and $b y^2=-1$.  If $-m-1$ is a non-square, then we are done; if $-m-1$ is a square then
we can continue the induction.  Hence we must also represent a non-square, and the proof is concluded.
\end{proof}

We are now ready to prove Theorems \ref{twoalmost}, \ref{onealmost}, and \ref{trialmost}.

\begin{proof}[Proof of Theorem \ref{twoalmost}]
Since $8T_x+1=(2x+1)^2$, the representation of $n$ by $f(x,y,z)=ax^2+by^2+cT_z$ is equivalent
to the representation of $8n+c$ by
$$
Q(x,y,z)=8ax^2+8by^2+cz^2
$$
with $z$ odd. The number of the latter representations
equals the number of solutions with $z$ arbitrary minus those with $z$ even.
As described in Section \ref{NotationSection}, every sufficiently large integer locally represented with bounded divisibility at the (finitely many) anisotropic primes of $Q'=Q$ or $Q'(x,y,2z)$ are represented, outside of the finitely many spinor exceptional square classes for $Q'(x,y,z)$ or $Q'(x,y,2z)$.  Thus, if the local conditions are satisfied, then
\begin{equation}\label{Xeqn}
\{ 8n+c: n\in  E(f)\} \subseteq \left(\bigcup_{j=1}^{r} \bigcup_{p\text{ anisotropic}}
\{ n_j p^{2s}: s \in \nature\}\right)  \cup  \left( \bigcup_{i=1}^{m} t_i\integer^2\right),
\end{equation}
where $p$ runs over the (finitely many) anisotropic primes, $n_1,\ldots,n_r$ are the finitely many ``sporadic'' natural numbers not
represented by $Q$, and $t_1\Z^2,\ldots,t_m\Z^2$ are finitely many spinor exceptional square classes which may not be represented.
Thus, $E(f)$ is a subset of a union of finitely many square classes, and hence
its asymptotic density is zero.

We then see that the local conditions at any odd prime $p$ are equivalent to those given in the theorem by
Lemma \ref{plocallemma} for $Q$.  We will use the original form $f$ to investigate the local condition at $p=2$.
A quick check shows that $T_z$ represents every integer modulo $8$ and Hensel's lemma then shows that $T_z$ represents
every integer $2$-adically.  Therefore, if $c$ is odd, then $cT_z$ represents every integer $2$-adically.
If $v_2(c)=1$, then $c T_z$ represents every even integer.  Since $\gcd(a,b,c)=1$, either $a$ or $b$ is odd,
hence every integer is $2$-adically represented.  If $v_2(c)=2$, then $cT_z$ represents every integer congruent to $4$
mod $8$.  Hence, we must represent $1$, $2$, and $3$ modulo $4$ with $ax^2+by^2$.  Without loss of generality,
we assume that $b$ is odd.  Then, $b$ is congruent to $1$ or $3$ modulo 4, and so is $by^2$ whenever $y$ is odd.
If $a$ is odd, then
either $a\equiv b\pmod{4}$, in which case $-b\pmod{4}$ is not represented, or $a\equiv -b\pmod{4}$,
in which case $2\pmod{4}$ is not represented.  Therefore, one sees $a\equiv 2\pmod{4}$,
which is equivalent to $v_2(a)=1$.  Finally, if $v_2(c)\geq 3$, then $ax^2+by^2$ cannot
represent every integer modulo 8 because an odd square is always congruent to 1 mod 8, so the local conditions are not satisfied.
\end{proof}

\medskip
\begin{proof}[Proof of Theorem \ref{onealmost}]
In this case the number of solutions to $n=f(x,y,z)=ax^2+bT_y+cT_z$
equals the number of representations of $8n+b+c$ by
$$
Q(x,y,z)=8ax^2+by^2+cz^2
$$
with $y$ and $z$ odd.  Thus, we again only need to show that every integer is locally represented.
The conditions given in the theorem for the odd primes are precisely those given by Lemma \ref{plocallemma}.
For $p=2$ we again use the fact that $T_y$ and $T_z$ represent every integer $2$-adically and note that
if $v_2(b)\leq 1$ or $v_2(c)\leq 1$ then every $2$-adic integer is represented because at least one of $a,b,c$
must be odd, and that if $2\leq v_2(c)\leq v_2(b)$ then not every integer is represented modulo 4.
\end{proof}

\medskip
\begin{proof}[Proof of Theorem \ref{trialmost}]
Clearly, $f(x,y,z)=aT_x+bT_y+cT_z$ represents the integer $n$ if and only if
$$
Q(x,y,z)=ax^2+by^2+cz^2
$$
represents $8n+a+b+c$ with $x,y,z$ all odd.  Again the local conditions at the odd primes are given
by Lemma \ref{plocallemma}.  For the $2$-adic conditions, we simply note that one of $a,b,c$ is odd, so every $2$-adic
integer is represented.
\end{proof}

\begin{proof}[Proof of Corollary \ref{quadreccor}]
We first note that if $b\R a$ then $\left(\frac{b}{a}\right)=1$. Thus, if the conditions given in
Theorems \ref{twoalmost}, \ref{onealmost}, or \ref{trialmost} hold,
then (by the multiplicative property of Jacobi symbols) we have
\begin{equation}\label{rsteqn}
1=\left(\frac{-b'c'}{a'}\right)\left(\frac{-a'c'}{b'}\right)\left(\frac{-a'b'}{c'}\right)
\left(\frac{2^r}{b'c'}\right)\left(\frac{2^s}{a'b'}\right)\left(\frac{2^t}{a'c'}\right),
\end{equation}
where $r,s,t$ are certain natural numbers.

By the law of quadratic reciprocity for Jacobi symbols,
\begin{align*}
&\left(\frac{-b'c'}{a'}\right)\left(\frac{-a'c'}{b'}\right)\left(\frac{-a'b'}{c'}\right)
\\=&\left(\frac{-1}{a'b'c'}\right)\cdot\left(\frac{b'}{a'}\right)\left(\frac{a'}{b'}\right)
\cdot \left(\frac{c'}{a'}\right)\left(\frac{a'}{c'}\right)\cdot\left(\frac{c'}{b'}\right)
\left(\frac{b'}{c'}\right)
\\=&(-1)^{\frac{a'-1}{2}+\frac{b'-1}{2}+\frac{c'-1}{2}}(-1)^{\frac{a'-1}{2}\cdot\frac{b'-1}{2}
+\frac{a'-1}{2}\cdot\frac{c'-1}{2}+\frac{b'-1}{2}\cdot\frac{c'-1}{2}}
\\=&(-1)^{\frac{a'+1}{2}\cdot\frac{b'+1}{2}\cdot \frac{c'+1}{2} -\frac{a'-1}{2}
\cdot\frac{b'-1}{2}\cdot \frac{c'-1}{2} - 1}
\end{align*}
Observe that $\frac{a'+1}{2}\cdot\frac{b'+1}{2}\cdot \frac{c'+1}{2}$ and
$\frac{a'-1}{2}\cdot\frac{b'-1}{2}\cdot \frac{c'-1}{2}$ have opposite parity if and only if
$a'\equiv b'\equiv c'\pmod{4}$.
So the product of three Jacobi symbols is $1$ if and only if $a'\equiv b'\equiv c'\pmod{4}$.

We finally deal with the $2$-power part.  If $a'\equiv b'\equiv -c'\pmod{8}$,
then $$1=\left(\frac{2}{a'b'}\right)=\left(\frac{2}{b'c'}\right)=\left(\frac{2}{a'c'}\right),$$
which concludes the first statement.  We now note that if $\pm a' \equiv c'+4\pmod{8}$,
then $\left(\frac{2}{a'}\right)=-\left(\frac{2}{c'}\right)$.
Therefore, in the cases $a'\equiv b'\equiv c'+4\pmod{8}$ or $\pm a'\equiv -b'\equiv c'+4\pmod{8}$ the Jacobi symbol
from the $2$-power part is $(-1)^{r+t}$, where $r$ and $t$ are as in equation (\ref{rsteqn}).
For (1), (4), and (5), we have $r=v_2(a)+1$, while $r=v_2(a)$ in the cases (2) and (3).
For (1), we have $t=v_2(b)+1$ and otherwise $t=v_2(b)$.  Thus for (1)-(3) we have $(-1)^{r+t}=(-1)^{v_2(a)+v_2(b)}$
and for (4)-(5) we have $(-1)^{r+t}=-(-1)^{v_2(a)+v_2(b)}$, from which we conclude the remaining two statements.
\end{proof}

\begin{proof}[Proof of Corollary \ref{new}] By Theorem \ref{twoalmost},
if $ax^2+by^2+2cT_z$ or $ax^2+cy^2+2bT_z$ is asymptotically universal then
we have (\ref{plocaleqn}).
Similarly,  if  $ax^2+2cy^2+bT_z$ or $ax^2+2by^2+cT_z$ is asymptotically universal then
we have
\begin{equation}\label{2} -bc\R a',\ -2ac\R b',\ -2ab\R c'.
\end{equation}

 Now assume that both $(\ref{plocaleqn})$ and $(\ref{2})$ hold. We want to deduce a contradiction.
 $(\ref{plocaleqn})$ and $(\ref{2})$ imply that $2\R b'$ and $2\R c'$. Recall that $v_2(b)\eq v_2(c)\ (\mo\ 2)$.
 So we have
\begin{equation}\label{3} -b'c'\R a',\ -a'c'\R b',\ -a'b'\R c'.
\end{equation}
It follows that
$$\left(\f{-b'c'}{a'}\right)=\left(\f{-a'c'}{b'}\right)=\left(\f{-a'b'}{c'}\right)=1.$$
Since $a'\eq b'\eq c'\ (\mo\ 4)$ fails, as in the proof of Corollary \ref{quadreccor} we have
$$\left(\f{-b'c'}{a'}\right)\left(\f{-a'c'}{b'}\right)\left(\f{-a'b'}{c'}\right)=-1.$$
So a contradiction occurs.
\end{proof}

  The following lemma gives a sufficient condition for a form not to be asymptotically universal.
It will be helpful for our proofs in the next section.
\begin{lemma}\label{vgt3lemma} For $f(x,y,z)=ax^2+by^2+cT_z, ax^2 + bT_y+cT_z, aT_x+bT_y+cT_z$ we define
$v_f=v_2(c),v_2(b+c),v_2(a+b+c)$ respectively.
If $v_f\geq 3$, then $f$ is not asymptotically universal.
\end{lemma}
\begin{proof}
Assume that $v_f\geq 3$ and $f$ is asymptotically universal. We want to deduce a contradiction.

In the case $f=ax^2+by^2+cT_z$, by Theorem \ref{twoalmost} we have $8\nmid c$ which contradicts $v_f\ge3$.

Now suppose that $f=ax^2 + bT_y+cT_z$. Since $4\nmid b$ or $4\nmid c$
by Theorem \ref{onealmost}, (up to symmetry) the vector $(b,c)$ modulo 8 is one of $(2,6)$, $(5,3)$, or $(1,7)$.
In the first case $a$ must be odd, while in the remaining two cases we have $bc\equiv 7\pmod{8}$
and hence $\left(\frac{2}{bc}\right)=1$.  Therefore, Theorem \ref{onealmost} and equation (\ref{rsteqn})
imply that $\left(\frac{-a'b'}{c'}\right)\left(\frac{-a'c'}{b'}\right)\left(\frac{-b'c'}{a'}\right)=1$.
However, the calculation from Corollary \ref{quadreccor} shows that
$$
\left(\frac{-a'b'}{c'}\right)\left(\frac{-a'c'}{b'}\right)\left(\frac{-b'c'}{a'}\right)=1
$$
if and only if $a'\equiv b'\equiv c'\pmod{4}$.  Since $b'\equiv 1\pmod{4}$ and $c'\equiv 3\pmod{4}$,
we are led to a contradction.

Finally we handle the case $f=aT_x+bT_y+cT_z$.
By Theorem \ref{trialmost} and $v_f\ge3$, the vector $(a,b,c)$ modulo 8 is one of $(8,1,7)$, $(8,3,5)$, $(2,5,1)$, $(2,3,3)$,
$(2,7,7)$, $(6,1,1)$, $(6,5,5)$, $(6,3,7)$, $(4,1,3)$, $(4,5,7)$.
The cases $(8,1,7)$ and $(8,3,5)$ are covered above.
 For the cases $(2,3,3)$, $(2,7,7)$, $(4,1,3)$, $(4,5,7)$, $(6,1,1)$, $(6,5,5)$ we have
 $$\left(\frac{-a'b'}{c'}\right)\left(\frac{-a'c'}{b'}\right)\left(\frac{-b'c'}{a'}\right)=-1\ \ \mbox{and}\ \
\left(\frac{2^{v_2(a)}}{b'c'}\right)=1,$$
 while in the cases $(2,5,1)$ and $(6,3,7)$
 we have
 $$\left(\frac{-a'b'}{c'}\right)\left(\frac{-a'c'}{b'}\right)\left(\frac{-b'c'}{a'}\right)=1
 \ \ \mbox{and}\ \ \left(\frac{2}{b'c'}\right)=-1.$$
 In view of (3.3), we get a contradiction.
\end{proof}

\section{On Almost Universal Forms}\label{SuffLargeSection}
In this section we investigate almost universal forms.  We will determine when asymptotically universal forms are
not almost universal.  We first consider sums with two squares.
\begin{proof}[Proof of Theorem \ref{twosufflarge}]
Assume the conditions of Theorem \ref{twoalmost}.  Recall that $n$ is represented by
$f(x,y,z)=ax^2+by^2+cT_z$ if and only $8n+c$ is represented by
$$
Q(x,y,z)=8ax^2+8by^2+cz^2
$$
with $z$ odd.  Since $v_2(c)\leq 2$  there are no representations of $8n+c$ by $Q(x,y,2z)$
due to congruence conditions modulo 8, thus the odd condition can be removed.  Therefore,
\begin{equation}\label{Xeqn2}
E(f)=\left\{ \frac{n-c}{8}:\ n\equiv c\ (\mo\ 8),\  Q(x,y,z)=n\text{ has no integral solution}\right\}.
\end{equation}
Let $t\integer^2$ be a spinor exceptional square class for the
genus of $Q$ such that $t$ is squarefree and $tx^2 \equiv
c\pmod{8}$ for some $x$.  We will see below that
$K=\rational(\sqrt{-tabc})$ will always be $\rational(i)$ or
$\rational(\sqrt{-2})$. Thus by the results of Earnest, Hsia, and
Hung \cite{EarnestHsia1} $t$ is a spinor exception for
the genus because $tx^2$ satisfying the Schulze-Pillot conditions
will imply that $t$ satisfies the Schulze-Pillot conditions.

When $t$ is not represented by the spinor genus of $Q$, Schulze-Pillot \cite{SchulzePillot2}
showed that for every prime $p$ splitting in $K$,
we have that $tp^2$ is not represented by the spinor genus of $Q$, and hence not by $Q$ (see \cite{SchulzePillot3}
for a full list of such properties).  If $t$ is represented by the spinor genus of $Q$, then for $p$ inert in $K$
we have that $tp^2$ is not primitively represented by the spinor genus of $Q$ \cite{SchulzePillot3}.
Here a primitive representation means that
$\gcd(x,y,z)=1$.  Thus, for squarefree $t$ represented by the spinor genus of $Q$ but not represented by $Q$,
$tp^2$ is not represented when $p$ is inert in $K$, as the number of representations of $tp^2$ equals the number of
(primitive) representations of $t$ plus the number of primitive representations of $tp^2$, and both of these are zero.
Hence, in either case we have seen that there are infinitely many integers in $t\Z^2$ not represented by $Q$
so that $E(f)$ is infinite if such a $t$ exists.

Next we show that if no such $t$ exists then $E(f)$ is finite.  By equations (\ref{Xeqn}) and (\ref{Xeqn2}),
if there is no such $t$ with $t x^2\equiv c \pmod{8}$ for some $x\in\Z$, then
$$
\{ 8n+c: n\in  E(f)\} \subseteq \left(\bigcup_{j=1}^{r} \bigcup_{p\text{ anisotropic}}
 \{ n_j p^{2s}: s \in \nature\}\right),
$$
where $n_1,\ldots,n_r$ are ``sporadic'' exceptions. Note that if every integer is represented $p$-adically
by the quadratic form $8ax^2+8by^2+cz^2$ then $p$ is not anisotropic.  Assume $p\mid c$ and fix an integer $n$.
Clearly, any $p$-adic solution to $Q(x,y,z)=n$ gives a solution to $Q(px,py,pz)=np^2$.
Since $Q$ satisfies the condition of Lemma \ref{plocallemma} and for any fixed $y\in\Z$ relatively prime to $p$
the equation $ax^2=np^2-by^2$ has a solution with $x$ relatively prime to $p$,
there are more solutions to the equation $Q(x,y,z)=np^2$ than to
the equation $Q(x,y,z)=n$, hence $p$ is not anisotropic.

Thus, the only possible anisotropic prime is $p=2$ and hence
$$
\{ 8n+c: n\in  E(f)\} \subseteq \bigcup_{j=1}^{r} \{ n_j 2^{2s}: s \in \nature\}.
$$
As $v_2(c)\leq 2$, we have
$$
\{ 8n+c: n\in  E(f)\} \subseteq \bigcup_{j=1}^{r} \{ n_j 2^{2s}: s \in \{0,1,2\} \},
$$
which shows that $E(f)$ is finite.

We now use Schulze-Pillot's classification \cite{SchulzePillot2} to determine the spinor exceptional square classes
$t_i\Z^2$.  Let a spinor exceptional square class $t\integer^2$ be given.  Earnest, Hsia, and Hung showed that if an odd prime $p$
is ramified in $K=\rational(\sqrt{-td})$ then $Q_p\iso \left\langle u_1,u_2p^r,u_3p^s\right\rangle$ with $u_i$
units in $\integer_p$ and $0<r<s$ (cf. \cite[Theorem 1(b)]{EarnestHsia1}).  However, since $p$ divides at most one of $a,b,c$,
this cannot occur.  It follows that $p$ is unramified in $K$, hence $K=K_{abc}$ or $K=K_{2abc}$.

Recall that $\SF(a'b'c')$ is the odd squarefree part of $abc$.
Assume that a prime $p$ dividing $\SF(a'b'c')$ is not split in
$K$.  Then, by Theorem 1(a) of Earnest, Hsia, and Hung
\cite{EarnestHsia1}, we have $Q_p\iso \left\langle u_1,u_2
p^{2r},u_3p^{2s}\right\rangle$ from the necessary condition given
by Kneser \cite{Kneser1}.  But this would contradict the fact that
$v_p(abc)$ is odd.  Conversely, when $p$ is split in $K$,
Earnest, Hsia and Hung showed that the local conditions for
$t$ to be a spinor exception are satisfied (cf.
\cite{EarnestHsia1}). If $p$ is odd and $v_p(abc)$ is even, then
\cite[Theorem 1(a)]{EarnestHsia1} shows that $t\not\eq0\ (\mbox{mod}\ p)$ satisfies
the necessary and sufficient conditions. Thus, the only possible
spinor exceptional square classes are given by $t=\SF(a'b'c')$ or
$t=2\SF(a'b'c')$. If $t\not\equiv 2^{-2s}c\pmod{8}$ for some
$s\in\nature$ with $2s\le v_2(c)$, then this spinor exceptional square class will not
occur in our consideration. Hence we conclude that $t=\SF(a'b'c)$.
Since $tabc$ times a suitable square equals $aa'bb'c^2$, we have
$K=K_{abc'}$.

From the above we see that every $p\mid \SF(a'b'c')$ must be split in $K$, which gives condition
(\ref{primeprodcond}).
If $Q$ represents $t$, then $Q$ also represents $t\integer^2$ (not necessarily primitively),
and hence condition (\ref{repcond}) is necessary.

We finally deal with the $2$-adic conditions.  Let $\beta$ be a
prime ideal of $O_K$ dividing $2$. Since $K=\rational(i)$ or
$K=\rational(\sqrt{-2})$, the $2$-adic completion
$K_{\beta}/\rational_2$ is ramified. After division by common
powers of $2$, we get
$$
Q_2\iso \left\langle c', 2^r b', 2^s a' \right\rangle,
$$
where $3-v_2(c)+v_2(b)=r\leq s=3-v_2(c)+v_2(a)$.

We now separate into cases depending on $v_2(c)$.  First we consider the case where $v_2(c)=2$.
In this case, we divide the equation $Q(x,y,z)=8n+c$ by 4 to find that
representation of $n$ by $f$ is equivalent to representation of $2n+1$ by $Q'(x,y,z)=2ax^2+2by^2+c' z^2$.
We recall by the conditions of Theorem \ref{twoalmost} that $v_2(a)=1$
and $v_2(b)=0$.  Thus, $L_2\iso \left\langle c',2b', 4 a'\right\rangle$.  Since $r\leq 3$ and $s\leq 2$,
the conditions of Theorem \ref{EarnestHsiaHungTheorem} are always satisfied for $v_2(4t)=2$.
Therefore every sufficiently large integer is always represented in this case.

For the remaining case $\SF(a'b'c)=t\equiv c\pmod{8}$, we conclude that $a'b'\equiv 1\pmod{2^{3-v_2(c)}}$.
Hence $a'\equiv b'\pmod{2^{3-v_2(c)}}$, which gives the first assertion of condition (\ref{abprimecond}).
Assume first that $c$ is odd.  When $r\geq 5$, Earnest and Hsia \cite{EarnestHsia2} showed that
$$
\theta(O^{+}(\left\langle a,b,c\right\rangle_2))=\theta(O^{+}(\left\langle ac,bc\right\rangle_2)).
$$
Since $a'\equiv b'\pmod{8}$ and scaling does not affect the spinor
norm, the lattice on the right-hand side is equivalent to
$\langle1,2^{s-r}\rangle_2$.  If $s=r$ then this is precisely
$N_2(\rational(i))$ as desired. Checking each case of Theorem
\ref{EarnestHsiaTheorem} shows that $\theta(O^{+}(\left\langle
1,2^{s-r} \right\rangle_2))\subseteq N_2(K)$, since $K=\rational(i)$
if $s-r$ is even, and $K=\rational(\sqrt{-2})$ if $s-r$ is odd.
Theorem \ref{EarnestHsiaHungTheorem} indicates that when $r\geq 5$
the sufficient conditions are also satisfied.
 Hence, if $4\mid b$ then $t$ is a spinor exception.  For $2\|b$ we have
 $5\in \theta(O^{+}(\left\langle c,2^{r}b'\right\rangle_2))\notin N_2(\rational(\sqrt{-2}))$;
it follows that $K=\rational(i)$ and hence $s$ is even.  But, when $r$ and $s$ have the same parity,
 none of the conditions of Theorem \ref{EarnestHsiaHungTheorem}(1) is satisfied when $r\geq 4$,
 therefore $t$ is a spinor exception.  For $r=3$, Theorem \ref{EarnestHsiaTheorem} implies that
 $K=\rational(\sqrt{-2})$ and $b\equiv c\pmod{8}$, and hence $s$ is even.
If $s=4$, then Theorem 2.2 of Earnest and Hsia \cite{EarnestHsia2} shows that
 $\theta(O^{+}(\left\langle a,b,c\right\rangle_2))=\rational_2\cross$.
 Therefore, $v_2(a)\geq 3$ is odd and $8\mid (b-c)$ in this case.  But then Theorem \ref{EarnestHsiaHungTheorem}(2)(c)
 is not satisfied since $s>2$, so it follows that $t$ is a spinor exception.

In the case $2\mid c$, we have $2\nmid b$. Thus we get
 $\langle c', 4b', 2^{s}a'\rangle_2$  after division by 2.
In view of the sublattice $\langle
c',4b'\rangle_2$, we have $K=\rational(i)$ and hence $2\mid s$. If
$s=2$ then taking the product of symmetries
$\sigma=\tau_{2x_1+x_2+x_3}\tau_{x_1}$ gives
$\theta(\sigma)=4(c'+b'+a')c'\notin N_2(\rational(i))$ because
each of $a',b',c'$ must be congruent to 1 mod $4$ by condition
(\ref{primeprodcond}).  Therefore $v_2(a)>0$ is even so that $2\mid
a$ and $v_2(a)\equiv c\pmod{2}$. In this case $L'=\langle
c',4b',2^s a'\rangle_2$ where $L'$ is as defined in Theorem
\ref{EarnestHsiaHungTheorem}(1), so $\theta(O^+(L')\subseteq
N_2(K)$. None of the conditions in Theorem
\ref{EarnestHsiaHungTheorem}(1)(c-d) can be satisfied, so $t$ is a
spinor exception.
\end{proof}

\begin{proof}[Proof of Corollary \ref{abc}]
 For any $n\in\N$ we have
\begin{align*}&2n+1=2ax^2+2by^2+cz^2\ \text{for some}\ x,y,z\in\Z
\\\iff&2n+1=2ax^2+2by^2+c(2z+1)^2\ \text{for some}\ x,y,z\in\Z
\\\iff&n-\frac{c-1}2=ax^2+by^2+4cT_z\ \text{for some}\ x,y,z\in\Z.
\end{align*}
By Theorem \ref{twosufflarge} and Theorem \ref{twoalmost},
\begin{align*}&ax^2+by^2+4cT_z\ \text{is almost universal}
\\\iff&ax^2+by^2+4cT_z\ \text{is asymptotically universal}
\\\iff&2\|ab,\ -8bc\ R\ a',\ -8ac\ R\ b',\ \text{and}\ -ab\ R\ c.
\end{align*}
So the first part of Corollary \ref{abc} follows.

When $n\in\N$, clearly
\begin{align*}&2n+1=2ax^2+c(y^2+z^2)\ \text{for some}\ x,y,z\in\Z
\\\iff&2n+1=2ax^2+4cy^2+cz^2\ \text{for some}\ x,y,z\in\Z.
\end{align*}
In the case $b=2c$,
\begin{align*}&2\|ab,\ -2bc\ R\ a',\ -2ac\ R\ b',\ \text{and}\ -ab\ R\ c
\\\iff&c=1,\ 2\nmid a,\ \text{and}\ -1\R a'.
\end{align*}
So we also have the second part of Corollary \ref{abc}.
\end{proof}

\begin{proof}[Proof of Corollary \ref{ab}]
 In light of Theorem \ref{twoalmost},
\begin{align*}& ax^2+by^2+2T_z\ \text{is asymptotically universal}
\\\iff&-b\R a'\ \text{and}\ -a\R b
\\\iff&ax^2+y^2+2bT_z\ \text{is asymptotically universal}
\end{align*}
and
\begin{align*}&ax^2+2y^2+bT_z\ \text{is asymptotically universal}
\\\iff&-b\R a'\ \text{and}\ -2a\R b
\\\iff&ax^2+2by^2+T_z \ \text{is asymptotically universal}
\end{align*}

 Now assume that $-a\R b$ and $-b\R a'$.
 Then both  $ax^2+by^2+2T_z$ and $ax^2+y^2+2bT_z$ are asymptotically universal.
 Recall that $\SF(a'b)=\SF(a')\SF(b)$ has a prime divisor $p\eq3\ (\mo\ 4)$.
 Whether $v_2(a)$ is even or odd, we cannot have
 both (1) and (2) of Theorem \ref{twosufflarge} for either of the two forms.
It follows that $ax^2+by^2+2T_z$ and $ax^2+y^2+2bT_z$ are almost universal.

Suppose that $-2a\R b$ and $-b\R a'$.
 Then both $ax^2+2y^2+bT_z$ and $ax^2+2by^2+T_z$ are asymptotically universal.
 As $2b\eq2\not\eq0\ (\mo\ 4)$ and not all prime divisors of $\SF(a'b)$ are congruent to 3 mod 4,
 we cannot have both (1) and (2) of Theorem \ref{twosufflarge} for either of the two forms.
So $ax^2+2y^2+bT_z$ and $ax^2+2by^2+T_z$ must be almost universal.
We are done.
\end{proof}

\begin{proof}[Proof of Corollary \ref{precisecor}]
(i) By Theorem \ref{twoalmost},
\begin{align*}& ax^2+y^2+T_z\ \text{is asymptotically universal}
\\\iff&-2\R a',\  \text{i.e.},\, \left(\frac{-2}p\right)=1\ \text{for each prime divisor}\ p\ \text{of}\ a'
\\\iff&\text{all odd prime divisors of }a\ \text{are congruent to 1 or 3 mod 8}.
\end{align*}
Similarly,
\begin{align*}& ax^2+2y^2+2T_z\ \text{is asymptotically universal}
\\\iff&\text{all prime divisors of }a\ \text{are congruent to 1 or 3 mod 8}.
\end{align*}

Now suppose that $-2\R a'$. As each prime $p\eq1,3\ (\mo\ 8)$ can be written in the form
$x^2+2y^2$ with $x,y\in\Z$, and
$$(x_1+2y_1^2)(x_2^2+2y_2^2)=(x_1x_2-2y_1y_2)^2+2(x_1y_2+x_2y_1)^2,$$
we can write $\SF(a')$ in the form $x_0^2+2y_0^2$ with $x_0,y_0\in\Z$ since
all prime divisors of $a'$ are congruent to 1 or 3 modulo 8.
If $a'\equiv 1\ (\mo\ 8)$, then $\SF(a')\eq1\ (\mo\ 8)$ and hence $y_0$ must be even,
so the equation $8(ax^2+y^2)+z^2=\SF(a')$ has a solution $(x,y,z)=(0,y_0/2,x_0)$,
which violates (3) in Theorem
\ref{twosufflarge} with $b=c=1$.
If $a'\eq3\ (\mo\ 8)$ then we don't have (1) in Theorem \ref{twosufflarge} with $b=c=1$.
Therefore, by Theorem
\ref{twosufflarge}, $ax^2+y^2+T_z$ must be almost universal.
When $a$ is odd, we have $4\nmid a$ and $v_2(2)\not\eq2\ (\mo\ 2)$, therefore
$ax^2+2y^2+2T_z$ is almost universal by Theorem \ref{twosufflarge} with $c=2$.

(ii) By Theorem \ref{twoalmost},
\begin{align*}& ax^2+2y^2+T_z\ (\text{or }ax^2+y^2+2T_z)\ \text{is asymptotically universal}
\\\iff&-1\R a',\  \text{i.e.},\, \left(\frac{-1}p\right)=1\ \text{for each prime divisor}\ p\ \text{of}\ a'
\\\iff&\text{all odd prime divisors of }a\ \text{are congruent to 1 mod 4}.
\end{align*}
Similarly,
\begin{align*}& ax^2+4y^2+2T_z\ (\text{or }ax^2+2y^2+4T_z)\ \text{is asymptotically universal}
\\\iff&\text{all prime divisors of }a\ \text{are congruent to 1 mod 4}.
\end{align*}

Below we assume that $-1\R a'$. It is well known that each prime $p\eq1\ (\mo\ 4)$
is a sum of two squares (of integers) and
$$(x_1+y_1^2)(x_2^2+y_2^2)=(x_1x_2-y_1y_2)^2+(x_1y_2+x_2y_1)^2.$$
So we can write $\SF(a')$ in the form $x_0^2+y_0^2$ with $x_0$ odd and $y_0$ even
(since all prime divisors of $a'$ are congruent to 1 mod 4). Thus the equation $4(ax^2+y^2)+z^2=\SF(a')$
has a solution $(x,y,z)=(0,y_0/2,x_0)$, which violates (3) in Theorem
\ref{twosufflarge} with $b=1$ and $c=2$. So $ax^2+y^2+2T_z$ is almost universal.
If $a'\equiv 1\ (\mo\ 8)$, then $\SF(a')\eq a'\not\eq5\ (\mo\ 8)$ and hence $4\mid y_0$,
so the equation $8(ax^2+2y^2)+z^2=\SF(a')$ has an integral solution $(x,y,z)=(0,y_0/4,x_0)$,
which violates (3) in Theorem
\ref{twosufflarge} for the form $ax^2+2y^2+T_z$.
If $a'\not\eq1\ (\mo\ 8)$ then we don't have (1) in Theorem \ref{twosufflarge} for the form $ax^2+2y^2+T_z$.
Thus, in view of Theorem
\ref{twosufflarge}, $ax^2+2y^2+T_z$ is also almost universal.

Now we also assume that $a$ is odd. Note that the equation $2(ax^2+2y^2)+z^2=\SF(a')$
has an integral solution $(x,y,z)=(0,y_0/2,x_0)$. Also, $v_2(4)\eq v_2(a)\ (\mo\ 2)$,
$$a\eq \SF(a)\eq1\ (\mo\ 8)\Longrightarrow 4(a0^2+4y^2)+z^2=\SF(a)\ \text{for some}\ y,z\in\Z,$$
and
$$a\eq\SF(a)\eq5\ (\mo\ 8)\Longrightarrow 4(ax^2+4y^2)+z^2=\SF(a)\ \text{for no}\ x,y,z\in\Z.$$
Thus, by Theorem \ref{twosufflarge},
 the form $ax^2+2y^2+4T_z$ is almost universal,
and $ax^2+4y^2+2T_z$ is almost universal if and only if $a\eq1\ (\mo\ 8)$.

The proof of is Corollary \ref{precisecor} is now complete.
\end{proof}

\begin{proof}[Proof of Corollary \ref{ExtendedCorollary}]
By Theorem \ref{twoalmost}, $ax^2+3y^2+T_z$ (resp., $ax^2+y^2+3T_z$, $ax^2+2y^2+6T_z$, $ax^2+6y^2+2T_z$)
is asymptotically universal if and only if both $-6\R a'$ and $a\eq1\ (\mo\ 3)$
(resp., $a\eq2\ (\mo\ 3)$, $a\eq1\ (\mo\ 6)$, $a\eq5\ (\mo\ 6)$).
Observe that $-6\R a'$ if and only if for each odd prime divisor $p$ of $a$ we have
$$\left(\f2p\right)=\left(\f{-3}p\right)=\left(\f p3\right),\ \text{i.e., }p\eq1,5,7,11\ (\mo\ 24).$$

For odd positive integers $b$ and $c$ not satisfying $a'\eq b\eq c\ (\mo\ 8)$, by Theorem \ref{twosufflarge},
the form $ax^2+bT_y+cT_z$ is almost universal if and only if it is asymptotically universal.
Thus $ax^2+3y^2+T_z$ (or $ax^2+y^2+3T_z$) is almost universal if and only if it is asymptotically universal.
If $a$ is odd then $4\nmid a$ and $v_2(2)=v_2(6)=1\not\eq6\eq2\ (\mo\ 2)$; thus by Theorem \ref{twosufflarge}
the form $ax^2+2T_y+6T_z$ (or $ax^2+6y^2+2T_z$) is almost universal if and only if it is asymptotically universal.

Combining the above, we have completed the proof of Corollary \ref{ExtendedCorollary}. \end{proof}

\begin{proof}[Proof of Corollary \ref{squarefree}]
Let $k,l\in\N$ with $k\ge l$. By Theorem \ref{twoalmost},
the form $2^kx^2+2^ly^2+mT_z$ is asymptotically universal if and only if $-2^{k+l}\R m'$ and
$$4\nmid m\ \ \text{or}\ \ (4\|m\ \&\ k=1\ \&\ l=0).$$

Assume that $2^kx^2+2^ly^2+mT_z$ is asymptotically universal.
As $v_2(m)<3$ and $2\nmid\SF(m')$, the equation
$$2^{3-v_2(m)}(2^kx^2+2^ly^2)+m'z^2=\SF(m')$$
has no integral solution if and only if $m'$ is not squarefree (i.e., $\SF(m')<m'$).
Thus, by Theorem \ref{twosufflarge},
the form $2^kx^2+2^ly^2+mT_z$ is not almost universal if and only if $k>0$ and $4\nmid m$ and
$$\begin{cases}l\le1\ \Longrightarrow\ k\eq m\ (\mo\ 2)&
\\l=0\ \&\ 2\nmid m\ \Longrightarrow\  k\ge3\ \&\ m\eq1\ (\mo\ 8).
\end{cases}$$

 In view of the above, we have the desired results in Corollary \ref{squarefree}.
\end{proof}

\begin{proof}[Proof of Corollary \ref{ExtendedCorollary2}]
By Theorem \ref{twoalmost},
\begin{align*}&ax^2+6^3y^2+T_z\ \text{is asymptotically universal}
\\\iff&-2^4 3^3 \R a'\ \text{and}\ -2a\R 3^3
\\\iff&-3\R a'\ \text{and}\ -2a\R 3
\\\iff &a\eq1\ (\mo\ 3)\ \text{and}\ \left(\f p3\right)=\left(\f{-3}p\right)=1
\ \text{for each prime divisor}\ p\ \text{of }a'
\\\iff&\text{all prime divisors of}\ a'\ \text{are congruent to 1 mod 3},\ \text{and}\ 2\mid v_2(a).
\end{align*}
And also,
\begin{align*}&ax^2+2\cdot 5^3y^2+T_z\ \text{is asymptotically universal}
\\\iff&-2^2 5^3 \R a'\ \text{and}\ -2a\R 5^3
\\\iff&-5\R a'\ \text{and}\ -2a\R 5
\\\iff &a\eq\pm 2\ (\mo\ 5)\ \text{and}\ \left(\f {-5}p\right)=1
\ \text{for each prime divisor}\ p\ \text{of }a'.
\end{align*}
For an odd prime $p$, clearly
$$\left(\f{-5}p\right)=1\iff \left(\f{-1}p\right)=\left(\f p5\right)
\iff p\eq 1,3,7,9\ (\mo\ 20)\iff 2\,\big|\,\left\lfloor\f p{10}\right\rfloor.$$

(i) Under the supposition, $ax^2+216y^2+T_z$ is asymptotically universal by the above.
If $8(ax^2+216y^2)+z^2=\SF(3^3a')=3\SF(a')$ for some $x,y,z\in\Z$, then we must have $x=0$ (since $8a>3a$)
and $3\mid z$, which contradicts that $3\nmid \SF(a')$. So the equation $8(ax^2+216y^2)+z^2=\SF(3^3a')$
has no integral solutions.
Applying Theorem \ref{twosufflarge} we find that $ax^2+216y^2+T_z$ is not almost universal if and only if
$a'\eq 3^3\ (\mo\ 8)$ and all prime divisors of $\SF(3^3a')=3\SF(a')$ are congruent to 1 or 3 mod 8.
Since $\SF(a')\eq a'\ (\mo\ 8)$ and
each prime divisor of $a'$ is congruent to 1 mod 3, the desired result follows.

(ii) Under the assumption, $a=2^{v_2(a)}a'\eq \pm2\ (\mo\ 5)$ and hence $ax^2+250y^2+T_z$ is asymptotically universal.
Note that the equation $8(ax^2+250y^2)+z^2=\SF(5^3a')=5\SF(a')$ has no integral solutions.
In view of Theorem \ref{twosufflarge}, $ax^2+250y^2+T_z$ is not almost universal if and only if
$a'\eq 5^3\eq5\ (\mo\ 8)$ and all prime divisors of $\SF(5^3a')=5\SF(a')$ are congruent to 1 mod 4.
Since $a'\eq\pm1\ (\mo\ 10)$ and
each prime divisor of $a'$ is congruent to one of $1,\, 3,\, 7,\, 9$ modulo 20, we finally obtain the desired result.
\end{proof}

\begin{proof}[Proof of Theorem \ref{onesufflargenec}]
As in Theorem \ref{twosufflarge}, $f$ will not be almost universal only if there is a (relevant)
anisotropic prime or a spinor exceptional square class with the correct congruence conditions modulo 8 for one of
the quadratic forms occuring in the inclusion/exclusion of theta series
$$\theta_{Q(x,y,z)}:=\theta_{Q'(x,y,z)}-\theta_{Q'(x,2y,z)}-\theta_{Q'(x,y,2z)}+\theta_{Q'(x,2y,2z)},$$
 where $Q'(x,y,z)=8ax^2+by^2+cz^2$.  We will first show that there are no relevant anisotropic primes.
 The conditions given by Theorem \ref{onealmost} imply that every odd prime $p$ is not anisotropic.
 By Lemma \ref{vgt3lemma}, the prime $2$ is never relevant because the congruence condition implies that
 the $2$-adic orders are at most 2.

Also as in Theorem \ref{twosufflarge}, the local conditions at each odd prime imply that the only
possible spinor exceptional square classes are $t\integer^2$ with $t=\SF(a'b'c')$ or $t=2\SF(a'b'c')$.
Moreover, the sufficient local conditions for the odd primes are satisfied if and only if
every prime divisor of $\SF(a'b'c')$ is split in $K=\rational(\sqrt{-2abc t})$.

If $t$ is a spinor exception for the genus of $Q'(x,y,z)$, then $t$ is
a spinor exception for the genus of $Q'(x,2y,z)$ and condition (\ref{onerepcond}) implies that $t$ is represented
the same number of times by each quadratic form.  If $t$ is not represented by the spinor genus of $Q'$ then $tp^2$
is also not (primitively) represented, where $p$ is an odd prime split in $K$.
If $t$ is represented by the spinor genus of $Q'$
then $tp^2$ is not primitively represented, where $p$ is an inert prime.  In either case $tp^2$ will also clearly not
be primitively represented by $Q'(x,2y,z)$, so $tp^2$ is not represented by $Q$.  Also,  if $t$ is not
a spinor exception for $Q'(x,y,z)$ or $Q'(x,2y,z)$ then $E(f)$ is finite.  Therefore,
for $E(f)$ to be infinite, it is sufficient that $t$ is
a spinor exception for the genus of $Q'$,
while it is necessary that $t$ is a spinor exception for the genus of $Q'(x,2y,z)$.

We now break into cases depending on $v_2(b+c)$.  Since $v_2(b+c)<3$ by Lemma \ref{vgt3lemma}, we begin
with the case $4\| b+c$.  Without loss of generality we assume $v_2(c)\leq v_2(b)$.
Since $v_2(c)<2$ by local conditions, congruence conditions imply in this case that $v_2(b)=v_2(c)\leq 1$.
If $v_2(c)=v_2(b)=1$, then $a$ is odd, and after division by common 2-powers we get $L=\langle c',b',4a'\rangle_2$,
so Theorem \ref{EarnestHsiaHungTheorem}(1) is not satisfied because $v_2(4t)=2\geq s$.
If $v_2(b)=v_2(c)=0$, then congruence conditions, without loss of generality,
give ($b\equiv 1\pmod 8\ \&\ c\equiv 3\pmod{8}$) or
($b\equiv 5\pmod{8}\ \&\ c\equiv 7\pmod{8}$).  In the second case we cannot have all prime divisors of $\SF(a'b'c')$
split in $K$.  In the first case, we must have $K=\rational(\sqrt{-2})$ and hence $s$ is odd.  But then $v_2(2^st)$ is
odd and Theorem \ref{EarnestHsiaHungTheorem}(2) implies that $0=v_2(b)>0$, thus $t$ is not a spinor exception.

For the remaining cases we note that $t\equiv b+c\pmod{8}$, so we must have $v_2(t)=v_2(b+c)$.
Conditions (\ref{onecongcond}), (\ref{oneprimeprodcond}), and (\ref{onerepcond}) now follow immediately.
First consider $v_2(b)\geq 5$ odd.  Then $\theta_Q=\theta_{Q'(x,y,z)}-\theta_{Q'(x,2y,z)}$,
where $Q'(x,y,z)=8ax^2+by^2+cz^2$.  If $c$ is even, then we have
$L_2=\left\langle c', 4a', 2^{r-1} b'\right\rangle_2$ and $K=\rational(i)$.  In this case, Earnest and Hsia proved
that the spinor norm can be considered only on $2\times 2$ sublattices.  Since $r-1$ is even,
Theorem \ref{EarnestHsiaTheorem} shows that $\theta(O^+(L_2))\subseteq N_2(K)$.  Moreover,
Theorem \ref{EarnestHsiaHungTheorem}(1)(c-d) cannot be satisfied, so $t$ is a spinor exception.
If $c$ is odd and $4\mid a$ then $r\geq 5$, so the spinor norm again equals the spinor norm on 2-dimensional sublattices.
Since $\SF(a'b'c')\equiv (b+c)'\equiv c'\pmod{8}$, we have $a'\equiv b'\pmod{8}$.  Therefore,
Theorem \ref{EarnestHsiaTheorem} shows that the spinor norm on each sublattice gives a subset of $N_2(K)$.
Since $r\geq 5$ none of the conditions of Theorem \ref{EarnestHsiaHungTheorem} is satisfied, and hence $t$ is a spinor
exception.  For $r=v_2(a)+3<5$, we have $\langle c', 2^r a', 2^s b'\rangle_2$.
Applying Theorem \ref{EarnestHsiaTheorem} to the the sublattice $\langle c', 2^r a'\rangle_2$,
we get $K=\rational(\sqrt{-2^r})$, and it follows that $s$ must be even.
 Therefore $t$ is not a spinor exception for $Q'(x,y,z)$ or $Q'(x,2y,z)$.

For $v_2(b)\geq 6$ even and $c$ even, in view of the sublattice
$\langle2c',8a\rangle_2$ and Theorem \ref{EarnestHsiaTheorem},
we have $K=\rational(i)$, which implies that $v_2(b)$ must be
odd.  For $v_2(b)\geq 6$ even and $c$ odd, Earnest and Hsia showed
that we may reduce the problem to 2-dimensional sublattices.  If
$a$ is odd, then the sublattice $\langle c,8a\rangle_2$ gives the
set $\{ \gamma: (\gamma,-2ac)_2=1\}$, which is a subgroup of
$N_2(\rational(\sqrt{-2}))$ if and only if $a\equiv c\pmod{8}$.
Theorem \ref{EarnestHsiaHungTheorem}(2)(c) shows that $t$ is a
spinor exception in this case, as $s>2$. When $a$ is even,  we
again note that $a'\equiv b'\pmod{8}$ by condition
(\ref{onecongcond}) and Theorem \ref{EarnestHsiaTheorem} implies
that $\theta(O^+(\langle8a,b,c\rangle_2))\subseteq N_2(K)$.

For $v_2(b)<3$, inclusion/exclusion gives $\theta_Q=\theta_{Q'}$.  For $b$ odd, the sublattice $\langle1,bc\rangle_2$
gives the spinor norm $Q(x_1+2x_2)Q(x_1){\rational_2\cross}^2=5{\rational\cross}^2\notin N_2(\rational(\sqrt{-2}))$,
 so $K=\rational(i)$ and $v_2(a)$ is even.  We then note that condition (\ref{onecongcond}) gives
 $1=\SF(a'b'c')\equiv (b+c)'\pmod{4}$, so that $b\equiv c\pmod{8}$.  For $s=3+v_2(a)\geq 5$,
the problem is now reduced to considering  2-dimensional sublattices, and we are done
since $\theta(O^{+}(\langle1,1\rangle_2))=N_2(\rational(i))$
  and $a'\equiv 1\pmod{4}$ by condition (\ref{oneprimeprodcond}).
  In this case $L''=L_2$ and $\theta(O^+(L'))\subseteq N_2(K)$, where $L'$ and $L''$ are as in
  Theorem \ref{EarnestHsiaHungTheorem}(1), so that condition (1)(d) is not satisfied and $t$ is a spinor exception.
  When $s=3$, Theorem \ref{EarnestHsiaTheorem} shows that
  $\theta(O^+(\langle8a,b,c\rangle_2))\not\subseteq N_2(\rational(i))$.
  For $v_2(b)=1$ Theorem \ref{EarnestHsiaTheorem} implies that $\theta(O^+(\langle2b',c'\rangle_2))\subseteq K$
  if and only if $K=\rational(\sqrt{-2})$ and $b'\equiv c'\pmod{8}$.  $K=\rational(\sqrt{-2})$ is equivalent to
  $2\nmid v_2(a)$ (which implies that $s\in\{4,6,\ldots\}$), while $b'\equiv c'\pmod{8}$ follows from
  $2b'+c\equiv \SF(a'b'c')\pmod{8}$ as each of $a',b',c'$ is congruent to $1$ or $3$ mod 8.
  For $s>4$ we are led to 2-dimensional
  lattices and Theorem \ref{EarnestHsiaTheorem} implies that
  $\theta(O^+(\langle8a,b,c\rangle_2))\subseteq N_2(\rational(\sqrt{-2}))$, while one sees that none of the conditions of
  Theorem \ref{EarnestHsiaHungTheorem}(2) can be satisfied, so $t$ is a spinor exception.
  For $s=4$, Theorem \ref{EarnestHsiaHungTheorem}(2)(a) is satisfied, so $t$ cannot be a spinor exception.

When $v_2(b)=2$, Theorem \ref{EarnestHsiaTheorem} implies that $K=\rational(i)$ and hence $v_2(a)$ is odd,
and Earnest and Hsia \cite{EarnestHsia2} showed that we may again consider the spinor norm on 2-dimensional sublattices,
to get $\theta(O^+(\langle8a,b,c\rangle_2))\subseteq N_2(K)$.
The conditions in Theorem \ref{EarnestHsiaHungTheorem}(1)(c-d) are not satisfied and $L'=L_2$, so $t$ is a spinor
exception.

For $v_2(b)=3$, if $t$ is a spinor exception for $Q'(x,y,z)$, then Theorem \ref{EarnestHsiaTheorem} for the sublattice $\langle1,8b'c'\rangle_2$
implies that $K=\rational(\sqrt{-2})$ and $b'\equiv c'\pmod{8}$.
Hence $v_2(a)$ must be odd.  If $v_2(a)=1$ then Theorem 2.2 of Earnest and Hsia \cite{EarnestHsia2} implies that
$\theta(O^+(\langle8a,b,c\rangle_2))=\rational_2\cross$.  For $v_2(a)>1$ odd we may again consider only
2-dimensional sublattices and Theorem \ref{EarnestHsiaHungTheorem}(2)(c) is satisfied, so $t$ is a spinor exception
for $Q'(x,y,z)$.  Finally, the property that $t$ is a spinor exception for $Q'(x,2y,z)$ is equivalent to the case where $r=5$
which was covered above.

For $v_2(b)=4$, if $t$ is a spinor exception for $Q'(x,y,z)$, then Theorem \ref{EarnestHsiaTheorem}
implies that $K=\rational(i)$ and hence $v_2(a)$ is odd.  Since $\SF(a'b'c')\equiv (b+c)'\equiv c'\pmod{8}$,
we have $a'\equiv b'\pmod{8}$ by condition (\ref{onecongcond})
and thus $\theta(O^+(\left\langle8a,b,c\right\rangle_2))\subseteq N_2(K)$.  Moreover, none of the conditions of
Theorem \ref{EarnestHsiaHungTheorem}(1) is satisfied, so $t$ is a spinor exception.  Finally,
the property that $t$ is a spinor exception
for $Q'(x,2y,z)$ is equivalent to the case where $r=6$, which was covered above.
\end{proof}

\begin{proof}[Proof of Corollary \ref{Tx+2Ty}] (i) By Theorem \ref{onealmost}, the form $ax^2+2T_y+T_z$
is asymptotically universal if and only if $-2\R\ a'$, i.e., each prime divisor of $a'$ is congruent to 1 or 3 mod 8.

Now assume that $-2\R a'$. As we mentioned before, $\SF(a')=2y^2+z^2$ for some $y,z\in\Z$. Clearly $z$ is odd.
If $\SF(a')\eq 2+1\ (\mo\ 2^3)$, then $y$ must be odd. Thus we cannot have both (1) and (3) in Theorem
\ref{onesufflargenec} with $b=2$ and $c=1$. Therefore $ax^2+2T_y+T_z$ is almost universal.

(ii) By Theorem \ref{onealmost}, the form $ax^2+4T_y+T_z$
is asymptotically universal if and only if $-1\R\ a'$, i.e., each prime divisor of $a'$ is congruent to 1 mod 4.

Now assume that $-1\R a'$. Then $\SF(a')=4y^2+z^2$ for some $y,z\in\Z$.
If $\SF(a')\eq 4+1\ (\mo\ 2^3)$, then $y$ must be odd. So we cannot have both (1) and (3) in Theorem
\ref{onesufflargenec} with $b=4$ and $c=1$. It follows that $ax^2+4T_y+T_z$ is almost universal. We are done.
\end{proof}

\begin{proof}[Proof of Corollary \ref{x^2+Ty}] Set $f_k(x,y,z)=2^k(x^2+T_y)+mT_z$ for $k=0,1,2,\ldots$.
By Theorem \ref{onealmost}, the form $f_k$
is asymptotically universal if and only if $-2\R m'$, and $2\nmid m$ when $k>0$.

 Assume that $f_k$ is asymptotically universal. Then all prime divisors of $m'$
 are congruent to 1 or 3 mod 8, thus $m'\eq1,3\ (\mo\ 8)$. Note that the equation
 $$8\times 2^kx^2+2^ky^2+mz^2=2^{v_2(2^k+m)}\SF(m')$$
 has no integral solutions with $yz$ odd since the right-hand side of the equation
 is smaller than $2^k+m$.

 {\it Case} 1. $k=0$. When $a=c=1$ and $b=m$, we obviously have $v_2(a)=0$, $2\nmid ac$ and $a\eq c\ (\mo\ 8)$.
 Thus, if $v_2(m)\not=4,6,\ldots$ then $f_0$ is almost universal by Theorem \ref{onesufflargenec}
 for the form $x^2+mT_y+T_z$.

 Now suppose $v_2(m)\in\{4,6,\ldots\}$. Then $\SF(m)\eq m+1\ (\mo\ 2)$. Also, $4\mid m+1$, $v:=v_2(m+1)=0$,
 $ m'\eq\SF(m')\ (\mo\ 8)$ and $(m+1)'=m+1\eq 1\ (\mo\ 8)$.
By Theorem \ref{onesufflargenec} for the form $x^2+mT_y+T_z$,
 $f_0$ is almost universal if $m'\eq3\ (\mo\ 8)$, and $f_0$ is not almost universal if $m'\eq1\ (\mo\ 8)$
 and $v_2(m)\not=4$.

 {\it Case} 2. $k>0$. In this case, $m$ is odd. If $k\in\{1,2\}$, then $f_k$ is almost universal since
 (4) in Theorem \ref{onesufflargenec} does not hold for $a=b=2^k$ and $c=m$. For $k\ge 3$, clearly $4\nmid 2^k+1$ and
 $\SF(m')=\SF(m)\eq m\eq 2^k+m=(2^k+m)'\ (\mo\ 8)$. Applying Theorem \ref{onesufflargenec}, we find that
 $f_k$ is not almost universal if $k>5$, or $k=3$ and $m\eq1\ (\mo\ 8)$.

 Combining the above, we have completed the proof.
\end{proof}

\begin{proof}[Proof of Corollary \ref{Tx+Ty}] Define $g_k(x,y,z)=2^k(x^2+2T_y)+mT_z$ for $k\in\N$.
By Theorem \ref{onealmost}, the form $g_k$
is asymptotically universal if and only if $-1\R m'$ (i.e., all prime divisors of $m'$ are congruent to 1 mod 4),
and $2\nmid m$ when $k>0$.

Suppose that $-1\R m'$. Clearly the equation
 $$8\times 2^kx^2+2^{k+1}y^2+mz^2=2^{v_2(2^{k+1}+m)}\SF(m')$$
 has no integral solutions with $y$ and $z$ odd (since the right-hand side of the equation
 is smaller than $2^{k+1}+m$). Note also that if $2\nmid m$ or $2\nmid v_2(m)$ then
 $\SF(2^k2^{k+1}m)=\SF(2m)\not\eq 2^{k+1}+m\ (\mo\ 2)$ and hence (2) in Theorem \ref{onesufflargenec}
 holds for the form $g_k$.

 {\it Case} 1. $k=0$. Since $v_2(1)-v_2(2)\not=2,4,6,\ldots$, if $v_2(m)\le v_2(2)=1$ (i.e., $4\nmid m$) then
 $x^2+2T_y+mT_z$ is almost universal by Theorem \ref{onesufflargenec}.
For $a=1$, $b=m$ and $c=2$, clearly $v_2(a)\not=1,3,5,\ldots$, and
$2\nmid a$ and $a\not\eq c\ (\mo\ 8)$. So, by Theorem
\ref{onesufflargenec} for the form $x^2+mT_y+2T_z$, $f_0$ is also
almost universal when $v_2(m)\in\{2,4,6,\ldots\}$. In the case
$v_2(m)\in\{3,5,\ldots\}$, clearly $\SF(m')\eq m'\eq1\eq
m/2+1=(m+2)'\ (\mo\ 4)$ and hence we have (1)-(4) in Theorem
\ref{onesufflargenec} for the form $x^2+mT_y+2T_z$. So Theorem
\ref{onesufflargenec} implies that $g_0$ is not almost universal if
$v_2(m)\in\{5,7,\ldots\}$.

{\it Case} 2. $k=1$. As $v_2(m+4)=0$ and $\SF(2'4'm')=\SF(m)\eq m\not\eq(m+4)'\ (\mo\ 2^3)$,
$g_1(x,y,z)=2x^2+4T_y+mT_z$ is almost universal by Theorem \ref{onesufflargenec}.

{\it Case} 3. $k\ge2$. In this case,
we have (1)-(3) in Theorem \ref{onesufflargenec} with $a=2^k$, $b=2^{k+1}$ and $c=m$.
Note also that $4\mid 2^k$. So, by Theorem \ref{onesufflargenec}, $g_k$ is not almost universal if $k>2$.

 In view of the above, we have proved both (i) and (ii) in Corollary \ref{Tx+Ty}.
\end{proof}

Now we turn to the proof of Theorem \ref{trisufflarge}.

\begin{proof}[Proof of Theorem \ref{trisufflarge}]
We assume without loss of generality that $v_2(a)\geq v_2(b)\geq v_2(c)=0$.
We again start by considering anisotropic primes, again arriving at the fact that only $p=2$ is possible.
However, Lemma \ref{vgt3lemma} implies bounded divisibility at $p=2$ by the congruence conditions, so there are
no relevant anisotropic primes.

We now determine when $t=\SF(a'b'c')$ or $t=2\SF(a'b'c')$ is a spinor exception.  For $E(f)$ to be infinite,
it is sufficient that $t$ is a spinor exception for $Q'(x,y,z)$, while it is necessary that $t$
is a spinor exception for one of the quadratic forms in the inclusion/exclusion.

We will break in to cases depending on $v:=v_2(a+b+c)$.  For $v=2$ and $v_2(a)<3$ we have $\theta_Q=\theta_{Q'}$
and $v_2(4t)=2\geq s=v_2(a)$ so none of the conditions of Theorem \ref{EarnestHsiaHungTheorem}
is satisfied and $t$ is not a spinor exception.  When $v=2$ and $v_2(a)\geq 3$,
we have $4\mid b+c$.  But then,  we may assume that $b\equiv 3\pmod{4}$ without loss of generality, as $b$ and $c$ are both odd.
Thus $K=\rational(\sqrt{-2})$ since every prime divisor of $b$ must split in $K$, and hence $v_2(tabc)$ is odd.
 However, Theorem \ref{EarnestHsiaHungTheorem}(2) implies that $v_2(tabc)$ must be even because $r=v_2(b)=0$.

We now must have $v\leq 1$, $t=2^v \SF(a'b'c')\equiv (a+b+c)\pmod{8}$, and $K=K_{2^v abc}$.
This gives conditions (\ref{tricongcond}), (\ref{triprimeprodcond}), and (\ref{trirepcond}).

For $v_2(b)\geq 5$ we have $\SF(a'b'c')\equiv (a+b+c)'\equiv c'\pmod{8}$, so $a'\equiv b'\pmod{8}$.
Again we are led to 2-dimensional sublattices and it follows that
$O^{+}(\langle a,b,c\rangle_2)\subseteq N_2(K)$, while none of the conditions of Theorem \ref{EarnestHsiaHungTheorem}
is satisfied, so $t$ is a spinor exception.

For $v_2(b)<3$ we have $\theta_{Q}=\theta_{Q'(x,y,z)}-\theta_{Q'(2x,y,z)}$.
For $b$ odd the sublattice $\langle b,c\rangle_2$ gives $b\equiv c\pmod{8}$ and $K=\rational(i)$,
so that $v_2(a)$ is odd.  But $b\equiv c\pmod{8}$ automatically by condition (\ref{tricongcond}).
If $s=v_2(a)\leq 3$ then Theorem 2.2 of Earnest and Hsia \cite{EarnestHsia2} implies that
$\theta(O^+(\langle a,b,c\rangle_2))=\rational_2\cross$.  so that $t$ is not a spinor exception for $Q'$.
For $s\geq 5$ Earnest and Hsia showed that we may reduce to 2-dimensional sublattices,
so that $\theta(O^+(\langle a,b,c\rangle_2))\subseteq N_2(K)$.  We then verify with
Theorem \ref{EarnestHsiaTheorem} that the Kneser condition is satisfied for $L'$ and $L''$ as defined in
Theorem \ref{EarnestHsiaHungTheorem}(1).  In this case condition 1(d) of Theorem \ref{EarnestHsiaHungTheorem}
is not satisfied, so $t$ is a spinor exception for $Q'$.  For $s=3$ and $x$ even, the situation
is similar to the case $s=5$, thus the above argument shows that $t$ is a spinor exception for $Q'(2x,y,z)$.

When $v_2(b)=1$, Theorem \ref{EarnestHsiaTheorem} for the sublattice $\langle c,2b'\rangle_2$ implies
that $K=\rational(\sqrt{-2})$
and $b'\equiv c'\pmod{8}$ (which is already satisfied by (\ref{tricongcond})).  Thus, $s$ is even.
If $s\leq 4$ then Theorem \ref{EarnestHsiaTheorem} for the sublattice $\langle c,2^s a'\rangle_2$
gives $5\in \theta(O^+(\langle a,b,c\rangle_2))$, so $t$ is not a spinor exception.
For $s>4$ we again split into 2-dimensional sublattices and the Kneser condition is satisfied by
Theorem \ref{EarnestHsiaTheorem}.  The Kneser condition for $L'$ as defined in Theorem \ref{EarnestHsiaHungTheorem}(2)
is satisfied and condition 2(c) is not satisfied, so $t$ is a spinor exception.
Again when $s=4$ we have $s'=6$ for $Q'(2x,y,z)$.

When $v_2(b)=2$, Theorem \ref{EarnestHsiaTheorem} for the sublattice $\langle c,4b'\rangle_2$ implies that
$K=\rational(i)$ and hence $s\geq 2$ is even.  For $s\geq 4$ we may reduce to 2-dimensional sublattices and the Kneser
condition is satisfied by Theorem \ref{EarnestHsiaTheorem}.
Theorem \ref{EarnestHsiaHungTheorem}(1)(c-d) are not satisfied, so $t$ is a spinor exception.
For $s=2$ we have $Q(2x_1 +x_2+x_3) Q(x_1)/4 \equiv 3 \pmod{4}$, so that we
don't have a spinor exception for $Q'$ in this case.  However, in the case
$s=2$, for $Q'(2x,y,z)$, it follows that $y$ is even by congruence
considerations and $r'=s'=4$.  As we will show later, $Q'(2x,y,z)$ has the
spinor exception $t$ in this case.

When $v_2(b)=3$, by Theorem \ref{EarnestHsiaTheorem} for the sublattice $\langle c,8b'\rangle_2$,
we have $K=\rational(\sqrt{-2})$
and hence $2\mid s$, as well as $b'\equiv c'\pmod{8}$.  For $s=4$, Theorem 2.2 of Earnest
and Hsia \cite{EarnestHsia2} implies that $\theta(O^+(L_2))=\rational_2\cross$, so $t$ is not a spinor exception.
When $s\geq 6$ we may again consider only 2-dimensional sublattices, and the Kneser condition is satisfied by
Theorem \ref{EarnestHsiaTheorem}.  Theorem \ref{EarnestHsiaHungTheorem}(2)(c) is not satisfied, so $t$ is a spinor
exception.  For $s=4$ and $x$ even we get a form with $s'=6$ and argue as above.  For $x, y$ even
we get $r'=v_2(b)+2\geq 5$ so that $t$ is a spinor exception.

Finally we deal with the case $v_2(b)=4$.  For $x, y$ even, we have $r'\geq 5$, so $t$ is a spinor exception
for $Q'(2x,2y,z)$.  Theorem \ref{EarnestHsiaTheorem} for the sublattice $\langle c,16b'\rangle_2$
implies that $K=\rational(i)$ and hence $2\mid s$.
But condition (\ref{tricongcond}) implies that $a'\equiv b'\pmod{8}$, so we find that for $s$ even,
$\theta(O^+(\langle a,b,c\rangle_2))\subseteq N_2(K)$.  Since $r=4$ and $s\geq 4$ is even, none of the conditions of
Theorem \ref{EarnestHsiaHungTheorem}(1) is satisfied, and $t$ is a spinor exception for $Q'(x,y,z)$.
\end{proof}

\begin{proof}[Proof of Corollary \ref{TTT}] By Theorem
\ref{trialmost}, the form $ax^2+2T_y+T_z$ is asymptotically
universal if and only if $-2\R a'$, i.e., each odd prime divisor
of $a$ is congruent to 1 or 3 modulo 8.

Below we assume that $-2\R a'$. As $\min\{v_2(a),v_2(2)\}\le1$,
(4) in Theorem \ref{trisufflarge} for the form $ax^2+2y^2+T_z$
just says that $v_2(a)-v_2(2)\in\{3,5,\ldots\}$.

Now suppose that $v_2(a)\in\{4,6,\ldots\}$. Then $v:=v_2(a+2+1)=0$
and
$$\SF(a'2'1')\eq(a+3)'\ (\mo\ 8)\iff a'\eq a+3\eq3\ (\mo\ 8).$$
Also, $\SF(2a)=2\SF(a')\not\eq a+2+1\ (\mo\ 2)$, and for any odd
integers $x,y,z$ we have $ax^2+2y^2+z^2>a\ge \SF(a'2'1')$.

 In view of the above, (1)-(4) in Theorem \ref{trisufflarge} for the form $ax^2+2y^2+T_z$
 are all valid if and only if $v_2(a)\in\{4,6,\ldots\}$ and $a'\eq
 3\ (\mo\ 8)$. Thus, by Theorem \ref{trisufflarge}, if $ax^2+2y^2+T_z$ is not almost universal,
 then we must have $a'\eq 3\ (\mo\ 8)$ and
 $v_2(a)\in\{4,6,\ldots\}$; when $v_2(a)\not=4$ (i.e.,
 $v_2(a)-v_2(2)\not=3$) the converse also holds.

 Combining the above we finally obtain the desired result.
\end{proof}

\end{document}